\documentclass[12pt,a4paper]{amsart}
\usepackage{t1enc}
\usepackage[latin1]{inputenc}
\usepackage[english]{babel}
\usepackage{hyperref}
\usepackage{graphicx}

\usepackage{latexsym}
\usepackage{amsmath,amsfonts,amssymb,amsthm}
\usepackage{mathptmx}

\begin{document}

\title{Volumes, Traces and Zeta Functions}
\author{Sergio Venturini}

\address{S. Venturini: Dipartimento Di Matematica, Universit\`{a} di Bologna, \,\,Piazza di Porta S. Donato 5 ---
I-40127 Bologna, Italy}
\email{venturin@dm.unibo.it}


\maketitle

%
%
\def\R{{\rm I\kern-.185em R}}
\def\C{{\rm\kern.37em\vrule height1.4ex width.05em depth-.011em\kern-.37em C}}
\def\N{{\rm I\kern-.185em N}}
\def\Z{{\bf Z}}
\def\Q{{\mathchoice{\hbox{\rm\kern.37em\vrule height1.4ex width.05em 
depth-.011em\kern-.37em Q}}{\hbox{\rm\kern.37em\vrule height1.4ex width.05em 
depth-.011em\kern-.37em Q}}{\hbox{\sevenrm\kern.37em\vrule height1.3ex 
width.05em depth-.02em\kern-.3em Q}}{\hbox{\sevenrm\kern.37em\vrule height1.3ex
width.05em depth-.02em\kern-.3em Q}}}}
\def\P{{\rm I\kern-.185em P}}
\def\H{{\rm I\kern-.185em H}}
%
\def\Aleph{\aleph_0}
\def\ALEPH#1{\aleph_{#1}}
\def\sset{\subset}\def\ssset{\sset\sset}
%
\def\bar#1{\overline{#1}}
\def\dim{\mathop{\rm dim}\nolimits}
\def\half{\textstyle{1\over2}}
\def\Half{\displaystyle{1\over2}}
\def\mlog{\mathop{\half\log}\nolimits}
\def\Mlog{\mathop{\Half\log}\nolimits}
\def\Det{\mathop{\rm Det}\nolimits}
\def\Hol{\mathop{\rm Hol}\nolimits}
\def\Aut{\mathop{\rm Aut}\nolimits}
\def\Re{\mathop{\rm Re}\nolimits}
\def\Im{\mathop{\rm Im}\nolimits}
\def\Ker{\mathop{\rm Ker}\nolimits}
\def\Fix{\mathop{\rm Fix}\nolimits}
\def\Res{\mathop{\rm Res}\nolimits}
\def\sp{\mathop{\rm sp}\nolimits}
\def\id{\mathop{\rm id}\nolimits}
\def\Trace{\mathop{\rm Tr}\nolimits}
\def\cancel#1#2{\ooalign{$\hfil#1/\hfil$\crcr$#1#2$}}
\def\prevoid{\mathrel{\scriptstyle\bigcirc}}
\def\void{\mathord{\mathpalette\cancel{\mathrel{\scriptstyle\bigcirc}}}}
\def\n{{}|{}\!{}|{}\!{}|{}}
\def\abs#1{\left|#1\right|}
\def\norm#1{\left|\!\left|#1\right|\!\right|}
\def\nnorm#1{\left|\!\left|\!\left|#1\right|\!\right|\!\right|}
%
\def\upperint{\int^{{\displaystyle{}^*}}}
\def\lowerint{\int_{{\displaystyle{}_*}}}
\def\Upperint#1#2{\int_{#1}^{{\displaystyle{}^*}#2}}
\def\Lowerint#1#2{\int_{{\displaystyle{}_*}#1}^{#2}}
%
\def\rem #1::#2\par{\medbreak\noindent{\bf #1}\ #2\medbreak}
\def\proclaim #1::#2\par{\removelastskip\medskip\goodbreak{\bf#1:}
\ {\sl#2}\medskip\goodbreak}
\def\ass#1{{\rm(\rmnum#1)}}
\def\assertion #1:{\Acapo\llap{$(\rmnum#1)$}$\,$}
\def\Assertion #1:{\Acapo\llap{(#1)$\,$}}
\def\acapo{\hfill\break\noindent}
\def\Acapo{\hfill\break\indent}
\def\proof{\removelastskip\par\medskip\goodbreak\noindent{\it Proof.\/\ }}
\def\prova{\removelastskip\par\medskip\goodbreak
\noindent{\it Dimostrazione.\/\ }}
\def\qed{{\bf //}\par\smallskip}
\def\BeginItalic#1{\removelastskip\par\medskip\goodbreak
\noindent{\it #1.\/\ }}
\def\iff{if, and only if,\ }
\def\sse{se, e solo se,\ }
\def\rmnum#1{\romannumeral#1{}}
\def\Rmnum#1{\uppercase\expandafter{\romannumeral#1}{}}
\def\smallfrac#1/#2{\leavevmode\kern.1em
\raise.5ex\hbox{\the\scriptfont0 #1}\kern-.1em
/\kern-.15em\lower.25ex\hbox{\the\scriptfont0 #2}}
%
\def\Left#1{\left#1\left.}
\def\Right#1{\right.^{\llap{\sevenrm
\phantom{*}}}_{\llap{\sevenrm\phantom{*}}}\right#1}
%
%
\def\a{\alpha}
\def\bg{\beta}
\def\g{\gamma}
\def\G{\Gamma}
\def\dg{\delta}
\def\D{\Delta}
\def\e{\varepsilon}
\def\eps{\epsilon}
\def\z{\zeta}
\def\th{\theta}
\def\T{\Theta}
\def\k{\kappa}
\def\lg{\lambda}
\def\Lg{\Lambda}
\def\m{\mu}
\def\n{\nu}
\def\r{\rho}
\def\s{\sigma}
\def\Sg{\Sigma}
\def\ph{\varphi}
\def\Ph{\Phi}
\def\x{\xi}
\def\om{\omega}
\def\Om{\Omega}


\newtheorem{theorem}{Theorem}[section]
\newtheorem{proposition}{Proposition}[section]
\newtheorem{lemma}{Lemma}[section]
\newtheorem{corollary}{Corollary}[section]
\newtheorem{remark}{Remark}[section]
\newtheorem{definition}{Definition}[section]

\newtheorem{teorema}{Teorema}[section]
\newtheorem{proposizione}{Proposizione}[section]
\newtheorem{corollario}{Corollario}[section]
\newtheorem{osservazione}{Osservazione}[section]
\newtheorem{definizione}{Definizione}[section]
\newtheorem{esempio}{Esempio}[section]
\newtheorem{esercizio}{Esercizio}[section]
\newtheorem{congettura}{Congettura}[section]

\bibliographystyle{abbrv}

\def\myZeta{\zeta}
\def\myD{D}
\def\myMatr{A}
\def\myMatrB{B}
\def\MatrSpace{M}
\def\MatrPSpace[#1]{{\MatrSpace^+}({#1},\R)}
\def\LyapMatrA{L}
\def\eigval{\lambda}
\def\mytrace{\alpha}
\def\mycbase{\beta}
\def\mycmatr{\gamma}
\def\mydecaya{\sigma}
\def\mydecayb{\tau}
\def\mydeg{a}
\def\myThetaSymb{\theta}
\def\myTheta[#1,#2,#3]{\myThetaSymb_{#1}\left({#2},{#3}\right)}
\def\myThetaStar[#1,#2,#3]{\myThetaSymb_{#1}^*\left({#2},{#3}\right)}
\def\myRes{\mathop{\rm Res}\limits}

\def\myXiSymb{\xi}
\def\myXi[#1,#2,#3]{\myXiSymb_{#1}\left({#2},{#3}\right)}
\def\myXiPlus[#1,#2,#3]{\myXiSymb_{#1}^+\left({#2},{#3}\right)}
\def\myXiMinus[#1,#2,#3]{\myXiSymb_{#1}^-\left({#2},{#3}\right)}

\def\myGbell{g}
\def\myvect#1{\underline{#1}}
\def\mya{\alpha}
\def\myamin{\alpha_{\min}}
\def\myamax{\alpha_{\max}}
\def\myMultiIndexA{p}
\def\myMultiIndexB{q}
\def\mydim{n}
\def\mydimB{p}
\def\myMinF{c_1}
\def\myMaxF{c_2}
\def\myConstA{c_3}
\def\myConstB{c_4}
\def\myindxA{i}
\def\myVarX{x}
\def\myVarY{y}
\def\myVarT{t}
\def\myVarTau{\tau}
\def\myIntVar{\om}
\def\myComplexVar{s}
\def\myHomVar{t}
\def\varRadius{r}
\def\nSmooth{m}
\def\nSmoothidx{k}
\def\myFuncA{\varphi}
\def\myFuncB{\psi}
\def\myDomainA{D}
\def\myDomain[#1]{\R^{#1}}
\def\myBall{B}
\def\myZMul{\lambda}
\def\myBdryBall{S}
\def\GenericSubset{E}
\def\FuncHolo{h}
\def\FuncReal{f}
\def\MyVector[#1,#2,#3]{({#1}_{#2},\ldots,{#1}_{#3})}
\def\myHomNorm[#1,#2]{\left|#2\right|_{#1}}
\def\myEucDot[#1,#2]{\left<{#1},{#2}\right>}
\def\Lyap{L}
\def\LyapName{Ljapunov}
\def\LyapBall{S_\Lyap}
\def\bigConstA{M}
\def\bigConstM{c}
\def\mygbell{g}

\def\myTranspose#1{{\left.{#1}\right.}^t}
\def\myDecaySpace[#1,#2]{{\mathcal{S}}_{#1}(\R^{#2})}
\def\myDecayNorm[#1,#2]{\norm{#1}_{#2}}
\def\mySymSobSpace[#1,#2,#3]{{\mathcal{S}}_{#1}^{#2}(\R^{#3})}

\def\myRadonA{\mu}
\def\myDirichlet[#1]{D_{#1}}

\def\myExpA{\lambda}
\def\myExpB{\mu}
\def\myE{e}

\def\myA{a}
\def\myB{b}
\def\myC{c}
\def\myK{k}
\def\myN{N}

\def\myNumA{\nu}
\def\rhVar{w}
\def\vLine{\e}
\def\vSector{\delta}
\def\vLineG{\sigma}

\begin{abstract}
Let $Q(x)$ be a quadratic form over $\mathbb{R}^n$.
The Epstein zeta function associated to $Q(x)$ is a 
well known function in number theory.
We generalize the construction of the Epstein zeta function to a class
of function $\varphi(x)$ defined in $\mathbb{R}^n$ that we call $A-$homogeneous,
where $A$ is a real aquare matrix of order $n$
having each eigenvalue in the left hal space $\Re\lambda>0$.
Such a class includes
all the homogeneous polynomials (positive outside the origin)
and all the norms on $\mathbb{R}^n$ which are smooth outside the origin.
As in the classical (i.e. quadratic) case we prove that such zeta functions
are obtained from the Mellin transforms of theta function of Jacobi type
associated to the $A-$homogeneous function $\varphi(x)$.
We prove that the zeta function associated to a $A-$homogeneous function $\varphi(x)$
which is positive and smooth outside the origin is an entire meromorphic function
having a unique simple pole at $s=\alpha$ the trace of the matrix $A$ with residue
given by the product of the trace $\alpha$ and the Lebesgue volume of the
unit ball associated to $\varphi(x)$, that is the volume of the set $x\in\R^n$
satisfying $\varphi(x)<1$.
We also prove that the theta funtion associated to $\varphi(x)$ has an asymptotic 
expansion near the origin. We find that the coefficients of such expansion
depend on the values that the zeta function associated to $\varphi(x)$ assumes
at the negative integers.

\end{abstract}

\section{\label{section:SectionIntro}Introduction}
For $\myComplexVar\in\C$ we denote by $\Re\myComplexVar$ and $\Im\myComplexVar$
respectively the real part and the imaginary part of $\myComplexVar$;
$\arg{\myComplexVar}$ is the principal determination of the
argument of $\myComplexVar$,
defined when $\myComplexVar$ is not a negative real number.

Let $\mydim$ be a positive integer,
$\myMatr$ a real square matrix of order $\mydim$,
and $\myHomVar>0$ a positive real number;
we put
\begin{equation}\nonumber
	\myHomVar^\myMatr:=e^{(\log t)\myMatr}
\end{equation}

We denote by
$\MatrPSpace[\mydim]$
the set of real square matrices of order $\mydim$
with $\Re\lambda>0$ for each engenvalue $\lambda$;
$I_\mydim$ will denote the identity matrix of order $\mydim$.

Given a fixed $\myMatr\in\MatrPSpace[\mydim]$
we say that a function $\myFuncA:\R^\mydim\to\R$
is $\myMatr-$\emph{homogeneous},
if for each $\myHomVar>0$
and each $\myVarX\in\R^\mydim$
\begin{equation}\label{def::HomogeneousFunction}
	\myFuncA(\myHomVar^\myMatr\myVarX)=\myHomVar\myFuncA(\myVarX).
\end{equation}

Then the \emph{unit ball} associated to $\myFuncA$ is
\begin{equation}\nonumber
	\myBall_\myFuncA=
		\bigl\{
			\myVarX\in\myDomain[\mydim]\mid\myFuncA(\myVarX)<1
		\bigr\}.
\end{equation}
For $\varRadius>0$ we also set
\begin{equation}\nonumber
	\myBall_\myFuncA(\varRadius)=
		\bigl\{
			\myVarX\in\myDomain[\mydim]\mid\myFuncA(\myVarX)<\varRadius
		\bigr\}.
\end{equation}

Let $\myMatr\in\MatrPSpace[\mydim]$ and
let $\myFuncA:\R^\mydim\to[0,+\infty[$
be an $\myMatr-$homogeneous function.
Setting $\myHomVar=2$ and $\myVarX=0$ in (\ref{def::HomogeneousFunction})
we obtain $\myFuncA(0)=2\myFuncA(0)$ and hence $\myFuncA(0)=0$.
We say that $\myFuncA$ is \emph{positive}
if $\myFuncA(\myVarX)>0$ when $\myVarX\neq0$.

Let us point out some examples of $\myMatr-$homogeneous functions.

Any homogeneous polynomial of degree $d$
is $d^{-1}I_\mydim-$homogeneous.

If $\myFuncA:\R^\mydim\to\R$ is 
the Minkowsky functional associated to an open star-like domain
$\myDomainA\sset\R^\mydim$ with respect to the origin
then $\myFuncA$ is the unique $I_\mydim-$homogeneous function
such that $\myBall_\myFuncA=\myDomainA$.

Let now $\myMatr\in\MatrPSpace[\mydim]$ and
let $\myFuncA:\R^\mydim\to[0,+\infty[$ be
a continuous positive $\myMatr-$homogeneous function.
Given a complex variable $s\in\C$ we define
\begin{equation}\label{eq::ZetaSeries}
	\myZeta(\myFuncA,\myComplexVar)=
	\sum_{\myIntVar\in\Z^\mydim\setminus\{0\}}\myFuncA(\myIntVar)^{-\myComplexVar},
\end{equation}
when the series on the right hand side converges.

We say that $\myZeta(\myFuncA,\myComplexVar)$ is the
$\myZeta$-function associated to $\myFuncA$.

The $\theta-$functions associated to $\myFuncA$ is defined
when $\tau$ is a complex number in the upper half plane $\Im\tau>0$
by the series
\begin{equation}\label{eq::ZetaSeries}
	\theta(\myFuncA,\tau)=
	\sum_{\myIntVar\in\Z^\mydim}e^{i\tau\myFuncA(\myIntVar)}
\end{equation}
and
\begin{equation}\label{eq::ZetaSeries}
	\theta^*(\myFuncA,\tau)=\theta(\myFuncA,\tau)-1=
	\sum_{\myIntVar\in\Z^\mydim\setminus\{0\}}e^{i\tau\myFuncA(\myIntVar)}.
\end{equation}

The purpose of this paper is to study
convergence and analytic continuation of such $\myZeta$-functions
and to give asymptotic expansion for the corresponding 
$\theta-$functions.

Let us fix now some further notations.

For each (finite) set $X$ we denote by $\#(X)$ the cardinality of $X$.

For each (measurable) subset $\GenericSubset\sset\R^\mydim$
we denote by $|\GenericSubset|$ the Lebesgue measure of $\GenericSubset$.

When $\myDomainA\sset\R^\mydim$ is an open domain
we denote (as usual) by $C^0(\myDomainA)$
(resp. $C^\nSmoothidx(\myDomainA)$ and $C^\infty(\myDomainA)$)
the space of the real continuous
(resp. differentiable of order $\nSmoothidx$
and indefinitely differentiable)
functions on $\myDomainA$.

The main results of this paper are the following:

\begin{theorem}\label{thm::Zetadef}
$\myMatr\in\MatrPSpace[\mydim]$ be a square matrix
and let $\mytrace$ be the trace of the matrix $\myMatr$.

Let $\myFuncA\in C^0(\R^\mydim)\cap C^{\infty}(\R^\mydim\setminus\{0\})$
be a continuous positive $\myMatr-$homogeneous function.
Then the series on the right hand side of (\ref{eq::ZetaSeries})
converges to a holomorphic function on the half space
$\Re\myComplexVar>\mytrace$ and extends to a holomorphic function
on $\C\setminus\{\mytrace\}$ having a simple pole at
$\myComplexVar=\mytrace$ with residue
\begin{equation}\nonumber
	\myRes_{\myComplexVar=\mytrace}=\mytrace|\myBall_\myFuncA|.
\end{equation}

We also have $\myZeta(\myFuncA,0)=-1$.
\end{theorem}

When $\mydim=1$ and $\myFuncA(x)=|x|$ then
$\myZeta(\myFuncA,\myComplexVar)=2\zeta(\myComplexVar)$,
where $\zeta(\myComplexVar)$ is the Riemann zeta function;
in this case the result is a classical one.

When $\mydim>1$ and $\myFuncA$ is a positive definite quadratic form then
$\myFuncA$ is $\myMatr-$homogeneous
for $\myMatr=2^{-1}I_\mydim$;
in this case $\myZeta(\myFuncA,\myComplexVar)$
is the Epstein zeta function associated to $\myFuncA$
and the result also is well known.

When $\myFuncA$ is a homogeneous polinomial of degree $d$
such that $\myFuncA(\myVarX)>0$ when $\myVarX\neq0$
then
$\myFuncA$ is $\myMatr-$homogeneous
for $\myMatr=d^{-1}I_\mydim$;
in this case
the meromorphic extension of $\myZeta(\myFuncA,\myComplexVar)$ 
has been established in
\cite{article:ChungMinAnBase}
but without an explicit computation of the residue
at the (unique) pole $\mydim/d$ of $\myZeta(\myFuncA,\myComplexVar)$.

\begin{theorem}\label{thm::Zetalim}
Let $\myMatr\in\MatrPSpace[\mydim]$ be a square matrix
with trace $\mytrace$.

Let $\myFuncA:\myDomain[\mydim]\to\R$ be a continuous
positive $\myMatr-$homogeneous function
Then:
\begin{enumerate}
	\item the series on the right hand side of (\ref{eq::ZetaSeries})
		converges when $\Re\myComplexVar>\mya$;
	\item if $\sigma\in\R$ we have
\begin{equation}\nonumber
	\lim_{\sigma\to\mya^+}
		(\sigma-\mya)\myZeta(\myFuncA,\sigma)
	=\mya|\myBall_\myFuncA|;
\end{equation}
	\item if $\varRadius>0$ we have
\begin{equation}\nonumber
	\lim_{\varRadius\to+\infty}
		\frac{\#(\myBall_\myFuncA(\varRadius)\cap\Z^\mydim)}{\varRadius^\mya}
	=|\myBall_\myFuncA|;
\end{equation}
\end{enumerate}
\end{theorem}

Our next general result on the $\myZeta$ functions
says that when $\myFuncA$ is not smooth throughout
$\R^\mydim\setminus\{0\}$ any result concerning the analytic continuation
of the serie $\myZeta(\myFuncA,\myComplexVar)$ cannot be given
simply by approximation with respect to the $C^0$ topology
(and so we are compelled to a tricky approach).

\begin{theorem}\label{thm::ZetaDiscontinuous}
Let $\myMatr\in\MatrPSpace[\mydim]$ be a square matrix
with trace $\mytrace$.

Let $\myFuncA:\myDomain[\mydim]\to\R$ be a continuous
positive $\myMatr-$homogeneous function
Given any $\dg>0$, $\e>0$ (arbitrarily small)
and $\bigConstA>0$ (arbitrarily large) there exists
$\myFuncB\in C^0(\R^\mydim)\cap C^{\infty}(\R^\mydim\setminus\{0\})$
such that for each $\myVarX\in\R^\mydim$
\begin{equation}\nonumber
	\myFuncA(\myVarX)\leq\myFuncB(\myVarX)\leq(1+\varepsilon)\myFuncA(\myVarX)
\end{equation}
but
\begin{equation}\nonumber
	\sup_{|\myComplexVar-\mytrace|=\delta}
	\left|\myZeta(\myFuncB,\myComplexVar)\right|\geq\bigConstA.
\end{equation}
\end{theorem}

The following theorem describes the behaviour of the $\theta-$function
associated to a continuous positive $\myMatr-$homogeneous function
smooth on $\R^\mydim\setminus\{0\}$.

If $\vLineG\in\R$ we adopt the following notation
for the Cauchy integral of a holomorphic function over a vertical line
\begin{equation}\nonumber
	\int_{(\vLineG)}f(\myComplexVar)d\myComplexVar
	=\int_{\vLineG-i\infty}^{\vLineG+i\infty}f(\myComplexVar)d\myComplexVar
	=\lim_{T\to+\infty}\int_{\vLineG-iT}^{\vLineG+iT}
		f(\myComplexVar)d\myComplexVar.
\end{equation}

\begin{theorem}\label{thm::ThetaAsymptotic}
Let $\myFuncA\in C^0(\R^\mydim)\cap C^{\infty}(\R^\mydim\setminus\{0\})$
be a positive $\myMatr-$homogeneous function.
Let $\rhVar$ be a complex number, and assume $\Re\rhVar>0$.
Let $\myN$ be a positive integer and let $0<\vLine<1$.
Then
\begin{eqnarray}\label{eq::ThetaAsymptotic}
\theta(\myFuncA,i\rhVar)&=&
\Gamma(\mytrace+1)\abs{\myBall_\myFuncA}\rhVar^{-\mytrace}+
\sum_{\myK=1}^{\myN}
	\frac{(-1)^\myK\myZeta(\myFuncA,-\myK)}{\myK!}\rhVar^\myK
\\
\nonumber
&&+
\frac{1}{2\pi i}\int_{(-\myN-1+\vLine)}\Gamma(\myComplexVar)\myZeta(\myFuncA,\myComplexVar)
	\rhVar^{-\myComplexVar}d\myComplexVar.
\end{eqnarray}
Moreover, given $0<\vSector<\pi/2$, there exists $\bigConstM>0$ which
depends on $\myN$, $\vLine$ and $\vSector$ only such that when
$\Re\rhVar>0$ and $\abs{\arg\rhVar}\leq\pi/2-\vSector$ then
\begin{equation}\label{eq::ThetaAsymptoticEstimate}
\abs{\int_{(-N-1+\vLine)}
	\Gamma(\myComplexVar)\myZeta(\myFuncA,\myComplexVar)
	\rhVar^{-\myComplexVar}d\myComplexVar}
\leq\bigConstM\abs{\rhVar}^{N+1-\vLine}.
\end{equation}
\end{theorem}

The paper is organized as follows.

In section \ref{section:SectionBasic} we fix some basic notations
and recall various known results that we need in the sequel of the paper.

In section \ref{section:SectionFuncEq} we study the Mellin transform
of a large class of theta function associated to any function $g(x)$
which satisfies $g(x)=O(\norm{x}^\mydecaya)$ as $\norm{x}\to+\infty$
and $\hat g(y)=O(\norm{y}^\mydecayb)$ as $\norm{y}\to+\infty$, where 
$\hat g(y)$ is the Fourier transform of $g(x)$ and $\mydecaya,\mydecayb>n$.

In section \ref{section:SectionHomog} we give a detailed description of
the behaviour of an $A-$homogeneous functions $\myFuncA(x)$
near the origin and for $\norm{x}\to\infty$.

In the remaining sections we give the proofs of the main results of this paper.

\section{\label{section:SectionBasic}Notations and some basic results}
We will denote by $\bigConstM_1,\bigConstM_2,\ldots$ suitable real positive constants.

If $f\in C^1(\R^\mydim)$ and $1\leq i\leq\mydim$ we denote by $\myD_if$ the derivative of
$f$ with respect to the variable $\myVarX_i$
and also set $\myD f=(\myD_1 f,\ldots,\myD_\mydim f).$

When $\myVarX=(\myVarX_1,\ldots,\myVarX_\mydim)\in\R^\mydim$
and $\myMultiIndexA=(\myMultiIndexA_1,\ldots,\myMultiIndexA_\mydim)\in\N^\mydim$
we set $\myVarX^\myMultiIndexA=
\myVarX_1^{\myMultiIndexA_1}\cdots\myVarX_\mydim^{\myMultiIndexA_\mydim}$
and $|\myMultiIndexA|=\myMultiIndexA_1+\cdots+\myMultiIndexA_\mydim$.
We also set
$\myD^\myMultiIndexA
=\myD_1^{\myMultiIndexA_1}\cdots\myD_\mydim^{\myMultiIndexA_\mydim}$,
as usual.

When $x,y\in\R^\mydim$ we denote by $\myEucDot[x,y]$
the Euclidean inner product of the vectors $x$ and $y$.

$GL(\mydim,\R)$ is the group of all invertible
real square matrices of order $\mydim$.

The transpose of the matrix $\myMatr$ is denoted by
$\myTranspose\myMatr$.

When $\myComplexVar$ is a complex variable $\Gamma(\myComplexVar)$ denotes the
Euler Gamma function.

Let us recall that $L^1(\R^\mydim)$ is the space of absolutely integrable
functions on $\R^\mydim$.

For $f\in L^1(\R^\mydim)$
\begin{equation}\nonumber
\hat{f}(y)=\int_\R f(x)e^{-2\pi i\myEucDot[x,y]}dx
\end{equation}
is the Fourier transform of $f$.

It is well known that if $\myMatr$ is any square real matrix
then, for some positive constant $\mycbase$,
\begin{equation}\nonumber	\norm{\myHomVar^\myMatr\myVarX}
\leq\bigConstM_1\myHomVar^{\mycbase}\norm{\myVarX}
\quad(\myVarX\in\R^\mydim,\ 1\leq\myHomVar<+\infty).
\end{equation}	

Let $\LyapMatrA$ be a positive definite symmetric matrix.
We say that $\LyapMatrA$ is a Ljapunov matrix
for the matrix $\myMatr$ if
\begin{equation}\label{eq::Lyap}
	\myEucDot[\myMatr\myVarX,\LyapMatrA\myVarX]>0
\end{equation}
for each $\myVarX\in\R^\mydim\setminus\{0\}$;

The following proposition give some characterizations
of the space $\MatrPSpace[\mydim]$.

\begin{proposition}\label{prop::PositiveMatrix}
Let $\mydim$ be a positive integer and let $\myMatr$
be a real square matrix of order $\mydim$.
Then the following conditions are equivalent.
\begin{enumerate}
	\item
		$\myMatr\in\MatrPSpace[\mydim]$;
	\item
		there exists a Ljapunov matrix for the matrix $\myMatr$;
	\item
		there exists a positive real constant $\mycmatr$ such that
		\begin{equation}\nonumber
			\myVarX\in\R^\mydim,\ 0<\myHomVar\leq1\ \Longrightarrow\ 
			\norm{\myHomVar^\myMatr\myVarX}\leq
			\bigConstM_1\myHomVar^{\mycmatr}\norm{\myVarX}.
		\end{equation}	
\end{enumerate}
\end{proposition}

For a proof of the proposition see e.g.
\cite{book:ArnoldOrdDiffEqSpringer}, section 22.

The proof of the following proposition is straightforward.

\begin{proposition}\label{prop::PBaseEstimates}
Let $\myMatr\in\MatrPSpace[\mydim]$.
Then there exist real constants $0<\mycmatr\leq\mycbase<+\infty$
such that for each $\myVarX\in\R^\mydim$
\begin{eqnarray}
	\label{eq::PBaseEstimatesOne}
	\myHomVar\geq1&\ \Longrightarrow\ &
		\bigConstM_1\myHomVar^\mycmatr\norm{\myVarX}
		\leq\norm{\myHomVar^\myMatr\myVarX}
		\leq\bigConstM_2\myHomVar^\mycbase\norm{\myVarX},\\
	\label{eq::PBaseEstimatesTwo}
	0<\myHomVar\leq1&\ \Longrightarrow\ &
		\bigConstM_1\myHomVar^\mycbase\norm{\myVarX}
		\leq\norm{\myHomVar^\myMatr\myVarX}
		\leq\bigConstM_2\myHomVar^\mycmatr\norm{\myVarX}.
\end{eqnarray}
\end{proposition}

The following proposition describes the ``$\myMatr-$polar coordinates''
when $\myMatr\in\MatrPSpace[\mydim]$.

\begin{proposition}\label{prop::PolarDecomposition}
Let $\myMatr\in\MatrPSpace[\mydim]$.
Let $0<\mycmatr\leq\mycbase<+\infty$ be as in Proposition \ref{prop::PBaseEstimates}.
Let $\Lyap$ be a
Ljapunov matrix for $\myMatr$
and set
$$
	\LyapBall=\bigl\{\myVarX\in\R^\mydim\vert\myEucDot[\Lyap\myVarX,\myVarX]=1\bigr\}.
$$
Then the map
$$
	S_\Lyap\times[0,+\infty[\ni(\bar\myVarX,\myHomVar)\mapsto
	\myVarX=\myHomVar^\myMatr\bar\myVarX\in\R^\mydim\setminus\{0\}
$$
is a diffeomorphism and when $\myVarX\in\R^\mydim\setminus\{0\}$,
$\bar\myVarX\in S_\Lyap$ and
$\myHomVar\in[0,+\infty[$ then
\begin{eqnarray}
\label{eq::PolarA}
\myHomVar\geq1,\ \myVarX=\myHomVar^\myMatr\bar\myVarX&\ \Longrightarrow\ &
\bigConstM_1\norm{\myVarX}^{\frac{1}{\mycbase}}\leq\myHomVar\leq
\bigConstM_2\norm{\myVarX}^{\frac{1}{\mycmatr}},\\
\label{eq::PolarB}
0<\myHomVar\leq1,\ \myVarX=\myHomVar^\myMatr\bar\myVarX&\ \Longrightarrow\ &
\bigConstM_3\norm{\myVarX}^{\frac{1}{\mycmatr}}\leq\myHomVar\leq
\bigConstM_4\norm{\myVarX}^{\frac{1}{\mycbase}}.
\end{eqnarray}
\end{proposition}

For a proof of the proposition see e.g.
\cite{book:ArnoldOrdDiffEqSpringer}, section 22.

Let $\myRadonA$ be a Radon measure on the half line $]0,+\infty[$
and let $\myComplexVar\in\C$, $a,b\in\R$, with $a<b$.

The Mellin transform of the measure $\myRadonA$ is defined formally by
the Mellin integral
\begin{equation}\nonumber
	\hat\myRadonA(\myComplexVar)=
	\int_{]0,+\infty[}\myHomVar^{\myComplexVar-1}d\myRadonA(\myHomVar).
\end{equation}

We also set
\begin{equation}\label{eq::MellinMinus}
	\hat\myRadonA^-(\myComplexVar)=
	\int_{]0,1[}\myHomVar^{\myComplexVar-1}d\myRadonA(\myHomVar),
\end{equation}
\begin{equation}\label{eq::MellinPlus}
	\hat\myRadonA^+(\myComplexVar)=
	\int_{[1,+\infty[}\myHomVar^{\myComplexVar-1}d\myRadonA(\myHomVar)
\end{equation}
and say that
$\hat\myRadonA^-(\myComplexVar)$ (resp. $\hat\myRadonA^+(\myComplexVar)$)
converges if the integral on the right hand side of
(\ref{eq::MellinMinus}) (resp. (\ref{eq::MellinPlus}))
is absolutely convergent.

When $a\in\R$
we say that the function
$\hat\myRadonA^-(\myComplexVar)$ (resp. $\hat\myRadonA^+(\myComplexVar)$) 
is defined on the half plane
$\Re\myComplexVar>a$ (resp. $\Re\myComplexVar<a$)
if there exists $a'\geq a$ (resp. $a'\leq a)$
such that the integral defining
$\hat\myRadonA^-(\myComplexVar)$ (resp. $\hat\myRadonA^+(\myComplexVar)$)
conveges absolutely when $\Re s>a'$ (resp. $\Re s<a')$
and defines a holomorphic function which extends
meromorphically on the larger half plane
$\Re\myComplexVar>a$ (resp. $\Re\myComplexVar<a$).

When $a,b\in\R$ satisfies $a<b$
we say that the Mellin transform $\hat\myRadonA(\myComplexVar)$
is defined on the strip $a<\Re\myComplexVar<b$ if
the function $\hat\myRadonA^-(\myComplexVar)$
is defined when $\Re\myComplexVar>a$,
and the function $\hat\myRadonA^+(\myComplexVar)$
is defined when $\Re\myComplexVar<b$;
in this case we define
\begin{equation}\nonumber
	\hat\myRadonA(\myComplexVar)=
	\hat\myRadonA^-(\myComplexVar)+
	\hat\myRadonA^+(\myComplexVar).
\end{equation}

We recall now Phragm\'en-Landau's theorem on Dirichlet integrals
(see e.g. \cite{book:TenenbaumAnalyticNumberTheory}, Theorem 6, pag. 111)

\begin{theorem}\label{thm::Landau}
Let $\myRadonA$ be a positive Radon measure on $\R$ having
support on $]\varepsilon,+\infty[$ for some $\varepsilon>0$.
Let $a,a'\in\R$, $a\leq a'$, and suppose that
\begin{equation}\nonumber
\myDirichlet[\myRadonA](\myComplexVar)=
\int_0^{+\infty}\myHomVar^{-\myComplexVar}d\myRadonA(\myHomVar)
\end{equation}
converges absolutely when $\Re\myComplexVar>a'$
and extends holomorphically on the larger half space
$\Re\myComplexVar>a$.
Then the integral $\myDirichlet[\myRadonA](\myComplexVar)$
also is absolutely convergent when $\Re\myComplexVar>a$.
\end{theorem}

Lastly we recall the Ikehara-Wiener theorem
(see e.g. \cite{book:LangAlgebraicNumerTheory}, pag. 305).

\begin{theorem}\label{thm::Ikehara}
Let $\myRadonA$ be a positive Radon measure on $\R$ having
support on $]\varepsilon,+\infty[$ for some $\varepsilon>0$
and for $\myVarX>0$ set
\begin{equation}\nonumber
	F(\myVarX)=\int_{[0,\myVarX[}d\myRadonA(\myHomVar)=\myRadonA([0,\myVarX[)
\end{equation}
Let $a\in\R$, $a\leq0$, and suppose that
\begin{equation}\nonumber
\myDirichlet[\myRadonA](\myComplexVar)=
\int_0^{+\infty}\myHomVar^{-\myComplexVar}d\myRadonA(\myHomVar)
\end{equation}
converges absolutely when $\Re\myComplexVar>a$
and extends to a meromorphic function in neighbourhood of the
region $\Re\myComplexVar\geq a$ having no pole
except for a simple pole at $\myComplexVar=a$
with residue $R$.
Then
\begin{equation}\nonumber
	\lim_{\myVarX\to+\infty}\frac{F(\myVarX)}{\myVarX^a}=
	\frac{R}{a}.
\end{equation}
\end{theorem}

\section{\label{section:SectionFuncEq}Functional equations}
{
This section is reminiscent of
the classical arguments related
to the functional equation of the Riemann zeta function.

Although in a totally different setting,
our approach looks formally like the Tate's thesis treatement
of the global functional equation for zeta function associated to
functions defined on adeles of a number field: 
compare Theorem 12 and Theorem 13 of e. g.
\cite{book:LangAlgebraicNumerTheory}, pag. 205-206)
with our theorem (\ref{thm::FunctionalEquation}).

From now to the end of the paper
$\mydim$ is a positive integer,
$\myMatr\in\MatrPSpace[\mydim]$,
$\mytrace>0$ is the trace of the matrix $\myMatr$ and
$0<\mycbase\leq\mycmatr<+\infty$ are constants satisfying
(\ref{eq::PBaseEstimatesOne}) and (\ref{eq::PBaseEstimatesTwo})
of proposition \ref{prop::PBaseEstimates}.

Let $\mydecaya$ and $\mydecayb$ two positive constants.

Given $\mygbell\in C^0(\R^\mydim)$ and $\myHomVar>0$ we put
\begin{eqnarray}
	\myTheta[\myMatr,\mygbell,i\myHomVar]&=&
	\sum_{\myIntVar\in\Z^\mydim}
	\mygbell\left(\myHomVar^\myMatr\myIntVar\right),
	\nonumber\\
	\myThetaStar[\myMatr,\mygbell,i\myHomVar]&=&
	\sum_{\myIntVar\in\Z^\mydim\setminus\{0\}}
	\mygbell\left(\myHomVar^\myMatr\myIntVar\right)=
	\myTheta[\myMatr,\mygbell,\myHomVar]-\mygbell(0)
	\nonumber
\end{eqnarray}
when the series on the right hand side converge absolutely.

Observe that such theta function are defined on
a (possibily empty) subset of
the half line of the complex plane
having zero real part and positive imaginary part.

When $\myThetaStar[\myMatr,\mygbell,i\myHomVar]$ is defined
for each $\myHomVar>0$ we denote the ``Mellin transform''
of the Radon measure
$\myThetaStar[\myMatr,\mygbell,i\myHomVar]d\myHomVar$
by
\begin{eqnarray}
	\myXi[\myMatr,\mygbell,\myComplexVar]&=&
	\int_0^{+\infty}\myThetaStar[\myMatr,\mygbell,i\myHomVar]
	\myHomVar^\myComplexVar\frac{d\myHomVar}{\myHomVar},
	\nonumber\\
	\myXiPlus[\myMatr,\mygbell,\myComplexVar]&=&
	\int_1^{+\infty}\myThetaStar[\myMatr,\mygbell,i\myHomVar]
	\myHomVar^\myComplexVar\frac{d\myHomVar}{\myHomVar},
	\nonumber\\
	\myXiMinus[\myMatr,\mygbell,\myComplexVar]&=&
	\int_0^1\myThetaStar[\myMatr,\mygbell,i\myHomVar]
	\myHomVar^\myComplexVar\frac{d\myHomVar}{\myHomVar},
	\nonumber
\end{eqnarray}
where $\myComplexVar$ is a complex variable
(see the discussion at the end of the previous section).

\begin{definition}
Let $\mydecaya>0$ be a positive real constant.
We denote by
\begin{equation}\nonumber
\myDecaySpace[\mydecaya,\mydim]
\end{equation}
the space of all continuous function
$\mygbell:\R^\mydim\to\C$ such that
\begin{eqnarray}\nonumber
&&\myDecayNorm[\mygbell,\mydecaya]=
	\sup_{\myVarX\in\R^\mydim}|\mygbell(\myVarX)|
	\left(1+\norm{\myVarX}^\mydecaya\right)
<+\infty,
\end{eqnarray}
\end{definition}

\begin{definition}
Let $\mydecaya>0$ and $\mydecayb>0$ be two positive real constants.
We denote by
\begin{equation}\nonumber
\mySymSobSpace[\mydecaya,\mydecayb,\mydim]
\end{equation}
the space of all continuous function
$\mygbell:\R^\mydim\to\C$ such that
$\myDecayNorm[\mygbell,\mydecaya]<+\infty$
and $\myDecayNorm[\hat\mygbell,\mydecayb]<+\infty$,
where $\hat\mygbell$ is the Fourier transform of $\mygbell$.
\end{definition}

Observe that $\mySymSobSpace[\mydecaya,\mydecayb,\mydim]$ endowed with
the norm $\norm{\mygbell}_{\mydecaya,\mydecayb}=
\myDecayNorm[\mygbell,\mydecaya]+\myDecayNorm[\hat\mygbell,\mydecayb]$
is a Banach space
and the Fourier transform is an isometry between the two spaces
$\mySymSobSpace[\mydecaya,\mydecayb,\mydim]$ and
$\mySymSobSpace[\mydecayb,\mydecaya,\mydim]$.
Moreover,
\begin{equation}\nonumber
\bigcap_{\mydecaya>0, \mydecayb>0}\mySymSobSpace[\mydecaya,\mydecayb,\mydim]
={\mathcal S}(\R^\mydim)
\end{equation}
is the usual Schwartz space of smooth function $\mygbell\in C^{\infty}(\R^\mydim)$
such that
\begin{equation}\nonumber
\sup_{\myVarX\in\R^\mydim}\abs{\myVarX^\myMultiIndexA D^\myMultiIndexB\mygbell(\myVarX)}<+\infty
\end{equation}
for each $\myMultiIndexA,\myMultiIndexB\in\N^\mydim$.

\begin{lemma}\label{prop::GBellEstimate}
Let $\mygbell\in\myDecaySpace[\mydecaya,\mydim]$.
Then, for each $\myHomVar>0$ and each $\myVarX\in\R^\mydim\setminus\{0\}$,
\begin{eqnarray}
\myHomVar\geq1\ &\Longrightarrow&\ 
\abs{\mygbell(\myHomVar^\myMatr\myVarX)}\leq
\bigConstM_1\myDecayNorm[\mygbell,\mydecaya]
\myHomVar^{-\mycmatr\mydecaya}\norm{\myVarX}^{-\mydecaya},
\nonumber\\
0<\myHomVar\leq1\ &\Longrightarrow&\ 
\abs{\mygbell(\myHomVar^\myMatr\myVarX)}\leq
\bigConstM_2\myDecayNorm[\mygbell,\mydecaya]
\myHomVar^{-\mycbase\mydecaya}\norm{\myVarX}^{-\mydecaya}.
\nonumber
\end{eqnarray}
where the constants $\bigConstM_1$ and $\bigConstM_2$
depend only on $\myMatr$.
\end{lemma}
\proof
Since $\mygbell\in\myDecaySpace[\mydecaya,\mydim]$ then
\begin{equation}\nonumber
\abs{\mygbell(\myHomVar^\myMatr\myVarX)}\leq
\frac{\myDecayNorm[\mygbell,\mydecaya]}{1+\norm{\myHomVar^\myMatr\myVarX}^\mydecaya}\leq
\frac{\myDecayNorm[\mygbell,\mydecaya]}{\norm{\myHomVar^\myMatr\myVarX}^\mydecaya}.
\end{equation}
By (\ref{eq::PBaseEstimatesOne}) and (\ref{eq::PBaseEstimatesTwo}) we have
$\norm{\myHomVar^\myMatr\myVarX}\geq\bigConstM_1\myHomVar^\mycmatr\norm\myVarX$
when $\myHomVar\geq1$ and 
$\norm{\myHomVar^\myMatr\myVarX}\geq\bigConstM_2\myHomVar^\mycbase\norm\myVarX$
when $0<\myHomVar\leq1$
and the assertion follows.
\qed

\begin{lemma}\label{prop::ThetaEstimates}
Let $\mygbell\in\myDecaySpace[\mydecaya,\mydim]$ and $\mydecaya>\mydim$.
Then the theta functions
$\myTheta[\myMatr,\mygbell,i\myHomVar]$ and
$\myThetaStar[\myMatr,\mygbell,i\myHomVar]$
are defined for each $\myHomVar>0$ and
\begin{eqnarray}
\label{eq::ThetaEstimateA}
\myHomVar\geq1\ &\Longrightarrow&\ 
\left|\myThetaStar[\myMatr,\mygbell,\myHomVar]\right|\leq
\bigConstM_1\myDecayNorm[\mygbell,\mydecaya]\myHomVar^{-\mycmatr\mydecaya},
\\
\label{eq::ThetaEstimateB}
0<\myHomVar\leq1\ &\Longrightarrow&\ 
\left|\myThetaStar[\myMatr,\mygbell,\myHomVar]\right|\leq
\bigConstM_2\myDecayNorm[\mygbell,\mydecaya]\myHomVar^{-\mycbase\mydecaya},
\end{eqnarray}
where the constant $\bigConstM_1$ and $\bigConstM_2$
depend only on $\myMatr$ and $\mydecaya$.
\end{lemma}

\proof
The assertion follows immediatly from the previous lemma, observing that
if $\mydecaya>\mydim$ then
\begin{equation}\nonumber
\sum_{\myIntVar\in\Z^\mydim\setminus\{0\}}\norm\myIntVar^{-\mydecaya}<+\infty.
\end{equation}
\qed

\begin{theorem}\label{thm::JacobiTransform}
Let $\mygbell\in\mySymSobSpace[\mydecaya,\mydecayb,\mydim]$
with $\mydecaya>\mydim$ and $\mydecayb>\mydim$.
Then the theta functions
$\myTheta[\myMatr,\mygbell,\myHomVar]$ and
$\myTheta[\myTranspose\myMatr,\hat\mygbell,\myHomVar]$
are defined for each $\myHomVar>0$ and satisfies the identity
\begin{equation}\label{eq::JacobiTransform}
\myTheta[\myMatr,\mygbell,\frac{i}{\myHomVar}]=
\myHomVar^\mytrace\myTheta[{\myTranspose\myMatr},\hat\mygbell,i\myHomVar].
\end{equation}
\end{theorem}

\proof
It suffices to aply the Poisson summation formula 
(see e.g. \cite{book:TenenbaumAnalyticNumberTheory}, Theorem 1, pag. 91,
or \cite{book:SteinAndWeiss}, Corollary 2.6, pag. 252)
to the function
$\mygbell_\myHomVar(\myVarX)=\mygbell\left(\myHomVar^{-\myMatr}\myVarX\right)$,
observing that
$\hat\mygbell_\myHomVar(\myVarX)=
	\myHomVar^\mytrace
	\hat\mygbell\bigl(\myHomVar^{\myTranspose\myMatr}\myVarX\bigl)$.
\qed

\begin{lemma}\label{prop::XiPlusEstimate}
Let $\mygbell\in\myDecaySpace[\mydecaya,\mydim]$ and $\mydecaya>\mydim$.

If $\Re\myComplexVar<\mycmatr\mydecaya$ then
$\myXiPlus[\myMatr,\mygbell,\myComplexVar]$ converges
and satisfies
\begin{equation}\label{eq::XiPlusEstimate}
\abs{\myXiPlus[\myMatr,\mygbell,\myComplexVar]}\leq
\frac{\bigConstM_1\myDecayNorm[\mygbell,\mydecaya]}{\mycmatr\mydecaya-\Re\myComplexVar}.
\end{equation}

If $\Re\myComplexVar>\mycbase\mydecaya$ then 
$\myXiMinus[\myMatr,\mygbell,\myComplexVar]$ converges
and satisfies
\begin{equation}\label{eq::XiMinusEstimate}
\abs{\myXiMinus[\myMatr,\mygbell,\myComplexVar]}\leq
\frac{\bigConstM_2\myDecayNorm[\mygbell,\mydecaya]}{\Re\myComplexVar-\mycbase\mydecaya}.
\end{equation}

The constant $\bigConstM_1$ and $\bigConstM_2$ depend only on $\myMatr$ and $\mydecaya$.
\end{lemma}

\proof
When $\Re\myComplexVar<\mycmatr\mydecaya$,
multiplying both sides of (\ref{eq::ThetaEstimateA})
by $\myHomVar^\myComplexVar$ and integrating, we obtain
\begin{eqnarray}
&&\abs{\myXiPlus[\myMatr,\mygbell,\myComplexVar]}\leq
	\int_1^{+\infty}\left|\myThetaStar[\myMatr,\mygbell,\myHomVar]
			\myHomVar^\myComplexVar\right|\frac{d\myHomVar}{\myHomVar}\leq
	\bigConstM_1\myDecayNorm[\mygbell,\mydecaya]
		\int_1^{+\infty}\myHomVar^{\Re\myComplexVar-\mycmatr\mydecaya}
			\frac{d\myHomVar}{\myHomVar}=
		\frac{\bigConstM_1\myDecayNorm[\mygbell,\mydecaya]}
			{\mycmatr\mydecaya-\Re\myComplexVar}
\nonumber
\end{eqnarray}
and (\ref{eq::XiPlusEstimate}) follows.
The proof of (\ref{eq::XiMinusEstimate}) is similar.
\qed

\begin{proposition}\label{prop::XiMinusEstimate}
Let $\mygbell\in\mySymSobSpace[\mydecaya,\mydecayb,\mydim]$
with $\mydecaya>\mydim$ and $\mydecayb>\mydim$.
Then $\myXiMinus[\myMatr,\mygbell,\myComplexVar]$
converges when $\Re\myComplexVar>\mytrace$ and
is defined when
$\Re\myComplexVar>\mytrace-\mycmatr\mydecayb$.
Moreover,
$\myXiPlus[\myTranspose\myMatr,\hat\mygbell,\mytrace-\myComplexVar]$
converges when $\Re\myComplexVar>\mytrace-\mycmatr\mydecayb$
and
\begin{equation}
\myXiMinus[\myMatr,\mygbell,\myComplexVar]=
	-\frac{\mygbell(0)}{\myComplexVar}
	-\frac{\hat\mygbell(0)}{\mytrace-\myComplexVar}
	+\myXiPlus[\myTranspose\myMatr,\hat\mygbell,\mytrace-\myComplexVar]
\end{equation}
\end{proposition}

\proof
After inserting $\myHomVar^{-1}$ in (\ref{eq::JacobiTransform})
we easily obtain
\begin{equation}\nonumber
\myThetaStar[\myMatr,\mygbell,i\myHomVar]=
-\mygbell(0)+\hat\mygbell(0)\myHomVar^{-\mytrace}+
\myHomVar^{-\mytrace}
\myThetaStar[\myTranspose\myMatr,\hat\mygbell,\frac{i}{\myHomVar}].
\end{equation}
By (\ref{eq::ThetaEstimateA}), when $0<\myHomVar\leq1$
\begin{equation}\nonumber
\abs{
	\myHomVar^{-\mytrace}
	\myThetaStar[\myTranspose\myMatr,\hat\mygbell,\frac{i}{\myHomVar}]
} \leq
\bigConstM_1\myDecayNorm[\hat\mygbell,\mydecayb]
	\myHomVar^{\mycmatr\mydecayb-\mytrace},
\end{equation}
and hence, when $\Re\myComplexVar>\mytrace$,
\begin{eqnarray}
\myXiMinus[\myMatr,\mygbell,\myComplexVar]&=&
-\mygbell(0)\int_0^1
	\myHomVar^\myComplexVar\frac{d\myHomVar}{\myHomVar}
+\hat\mygbell(0)\int_0^1
	\myHomVar^{\myComplexVar-\mytrace}\frac{d\myHomVar}{\myHomVar}
+\int_0^1\myThetaStar[\myTranspose\myMatr,\hat\mygbell,\frac{i}{\myHomVar}]
	\myHomVar^{\myComplexVar-\mytrace}\frac{d\myHomVar}{\myHomVar}
\nonumber\\
	&=&
	-\frac{\mygbell(0)}{\myComplexVar}
	-\frac{\hat\mygbell(0)}{\mytrace-\myComplexVar}
+\int_0^1\myThetaStar[\myTranspose\myMatr,\hat\mygbell,\frac{i}{\myHomVar}]
	\myHomVar^{\myComplexVar-\mytrace}\frac{d\myHomVar}{\myHomVar},
\nonumber
\end{eqnarray}
where all the integrals are absolutely convergent
if $\Re\myComplexVar>\mytrace$.

Making the change of variable $\myHomVar\to1/\myHomVar$
in the last integral we obtain
\begin{equation}\nonumber
\int_0^1\myThetaStar[\myTranspose\myMatr,\hat\mygbell,\frac{i}{\myHomVar}]
	\myHomVar^{\myComplexVar-\mytrace}\frac{d\myHomVar}{\myHomVar}=
\int_1^{+\infty}\myThetaStar[\myTranspose\myMatr,\hat\mygbell,i\myHomVar]
	\myHomVar^{\mytrace-\myComplexVar}\frac{d\myHomVar}{\myHomVar}
=\myXiPlus[\myTranspose\myMatr,\hat\mygbell,\mytrace-\myComplexVar].
\end{equation}
To end the proof it is enough to note that, by the previous lemma,
$\myXiPlus[\myTranspose\myMatr,\hat\mygbell,\mytrace-\myComplexVar]$
converges when
$\Re\myComplexVar>\mytrace-\mycmatr\mydecayb$.
\qed

\begin{theorem}\label{thm::FunctionalEquation}
Let $\mygbell\in\mySymSobSpace[\mydecaya,\mydecayb,\mydim]$
with $\mydecaya>\mydim$ and $\mydecayb>\mydim$.
Assume also that
\begin{equation}\label{eq::StripCondition}
\mydecaya+\mydecayb>\frac{\mytrace}{\mycmatr}.
\end{equation}
Then the  functions
$\myXi[\myMatr,\mygbell,\myComplexVar]$
and
$\myXi[\myTranspose\myMatr,\hat\mygbell,\myComplexVar]$
are defined in the strip
\begin{equation}\label{eq::StripDefinition}
\mytrace-\mycmatr\mydecayb<\Re\myComplexVar<\mycmatr\mydecaya
\end{equation}
and satisfy the identities
\begin{equation}\label{eq::IntegralEquation}
\myXi[\myMatr,\mygbell,\myComplexVar]=
-\frac{\mygbell(0)}{\myComplexVar}
-\frac{\hat\mygbell(0)}{\mytrace-\myComplexVar}
+\myXiPlus[\myMatr,\mygbell,\myComplexVar]
+\myXiPlus[\myTranspose\myMatr,\hat\mygbell,\mytrace-\myComplexVar]
\end{equation}
and
\begin{equation}\label{eq::FunctionalEquation}
\myXi[\myMatr,\mygbell,\myComplexVar]=
\myXi[\myTranspose\myMatr,\hat\mygbell,\mytrace-\myComplexVar].
\end{equation}
\end{theorem}

\proof
Condition (\ref{eq::StripCondition}) ensures that the strip
defined in (\ref{eq::StripDefinition}) is not empty.

By proposition (\ref{prop::XiMinusEstimate}),
$\myXiMinus[\myMatr,\mygbell,\myComplexVar]$ is defined
when $\Re\myComplexVar>\mytrace-\mycmatr\mydecayb$
and satisfies
\begin{equation}\nonumber
\myXiMinus[\myMatr,\mygbell,\myComplexVar]=
	-\frac{\mygbell(0)}{\myComplexVar}
	-\frac{\hat\mygbell(0)}{\mytrace-\myComplexVar}
	+\myXiPlus[\myTranspose\myMatr,\hat\mygbell,\mytrace-\myComplexVar]
\end{equation}

Adding $\myXiPlus[\myMatr,\mygbell,\myComplexVar]$,
which by lemma (\ref{prop::XiPlusEstimate}) converges
when $\Re\myComplexVar<\mycmatr\mydecaya$,
we obtain immediatly (\ref{eq::IntegralEquation}).

Replacing $\mygbell$ with its Fourier transform $\hat\mygbell$,
$\myMatr$ with its transpose and
$\myComplexVar$ with $\mytrace-\myComplexVar$
in equation (\ref{eq::IntegralEquation}) we easily obtain
\begin{equation}\label{eq::IntegralEquationDual}
\myXi[\myTranspose\myMatr,\hat\mygbell,\mytrace-\myComplexVar]=
-\frac{\tilde\mygbell(0)}{\myComplexVar}
-\frac{\hat\mygbell(0)}{\mytrace-\myComplexVar}
+\myXiPlus[\myMatr,\tilde\mygbell,\myComplexVar]
+\myXiPlus[\myTranspose\myMatr,\hat\mygbell,\mytrace-\myComplexVar],
\end{equation}
where $\tilde\mygbell$ is the function defined for each
$\myVarX\in\R^\mydim$ as
$\tilde\mygbell(\myVarX)=\mygbell(-\myVarX)$.

Of course $\tilde\mygbell(0)=\mygbell(0)$, and
by the symmetry with respect to the origin of
$\Z^\mydim\setminus\{0\}$ obviously
\begin{equation}\nonumber
\myThetaStar[\myMatr,\tilde\mygbell,i\myHomVar]
=\myThetaStar[\myMatr,\mygbell,i\myHomVar]
\end{equation}
which implies
\begin{equation}\nonumber
\myXiPlus[\myMatr,\tilde\mygbell,\myComplexVar]
=\myXiPlus[\myMatr,\mygbell,\myComplexVar]
\end{equation}
and so obtaining
\begin{equation}\label{eq::IntegralEquationBis}
\myXi[\myTranspose\myMatr,\hat\mygbell,\mytrace-\myComplexVar]=
-\frac{\mygbell(0)}{\myComplexVar}
-\frac{\hat\mygbell(0)}{\mytrace-\myComplexVar}
+\myXiPlus[\myMatr,\mygbell,\myComplexVar]
+\myXiPlus[\myTranspose\myMatr,\hat\mygbell,\mytrace-\myComplexVar].
\end{equation}

The identity (\ref{eq::FunctionalEquation}) 
follows now by comparing 
(\ref{eq::IntegralEquation}) and (\ref{eq::IntegralEquationBis}).

\qed

}
\section{\label{section:SectionHomog}$\myMatr$-Homogeneous functions}
{
\def\myDomain{\Omega}
\def\myF{f}
\def\myG{g}
\def\myM{\Lambda}
\def\myME{\lambda}
\def\myk{k}
\def\myi{i}
\def\myj{j}
\def\myl{l}

The purpose of this section is to prove that if
$\myFuncA\in C^0(\R^\mydim)\cap C^{\infty}(\R^\mydim\setminus\{0\})$
is $\myMatr-$homogeneous and positive then the functions of the form
$e^{-\myFuncA^\myExpA}$ and $\myFuncA^\myExpB e^{-\myFuncA}$
are in $\mySymSobSpace[\mydecaya,\mydecayb,\mydim]$
for each $\mydecaya>0$ and for suitable $\mydecayb>\mydim$
when the exponents $\myExpA$ and $\myExpB$ are large enought.

From elementary calculus we have the following:

\begin{lemma}\label{lemma::LinearChange}
Let $\myDomain\in\R^\mydim$ be a (open) domain and let
$\myF\in C^1(\myDomain)$.
Let $\myM=(\myME_{\myi\myj})\in GL(\mydim,\R)$
and set $\myG=\myF\circ\myM^{-1}$.
Then, for each $\myVarX\in\myDomain$,
\begin{equation}\label{lemma::LinearChange::MatrForm}
\myD\myF(\myM^{-1}\myVarX)=D\myG(\myVarX)\myM,
\end{equation}
that is, for $\myj=1,\ldots,\mydim$,
\begin{equation}\label{lemma::LinearChange::Extended}
\myD_\myj\myF(\myM^{-1}\myVarX)=
\sum_{\myi=1}^\mydim\myD_\myi\myG(\myVarX)\myME_{\myi\myj}.
\end{equation}
\end{lemma}


\begin{proposition}\label{prop::EulerEx}
Let $\myME\in\R$ and $\myM\in GL(\mydim,\R)$.
Let $\myk>0$ be an integer and let
$\myF\in C^\myk(\R^\mydim\setminus\{0\})$.
Suppose that for each $\myVarX\in\R^\mydim\setminus\{0\}$
\begin{equation}\label{prop::EulerEx::Homo}
	\myME\myF(\myVarX)=\myF(\myM^{-1}\myVarX).
\end{equation}
Then, given $\myk$ integers
$1\leq\myj_1,\ldots,\myj_\myk\leq\mydim$,
\begin{equation}\label{prop::EulerEx::ToProve}
	\myD_{\myj_1}\cdots\myD_{\myj_\myk}\myF(\myM^{-1}\myVarX)
	=\sum_{1\leq\myi_1,\ldots,\myi_\myk\leq\mydim}
	\myD_{\myi_1}\cdots\myD_{\myi_\myk}\myF(\myVarX)\myME
	\prod_{\myl=1}^\myk\myME_{\myi_\myl\myj_\myl}.
\end{equation}
\end{proposition}

\proof
The proof is by induction on $\myk$.
Using (\ref{prop::EulerEx::Homo}) and
(\ref{lemma::LinearChange::Extended}),
\begin{eqnarray}
\nonumber
\myD_{\myj_1}\myF(\myM^{-1}\myVarX)&=&
	\sum_{\myi_1=1}^\mydim
	\myD_{\myi_1}(\myF\circ\myM^{-1})(\myVarX)\lambda_{\myi_1\myj_1}
\\
\nonumber
&=&\sum_{\myi_1=1}^\mydim\myD_{\myi_1}\myF(\myVarX)\myME\myME_{\myi_1\myj_1},
\end{eqnarray}
which is just (\ref{prop::EulerEx::ToProve}) when $\myk=1$.

Assume that (\ref{prop::EulerEx::ToProve}) holds for $\myk-1$,
that is
\begin{equation}\nonumber
	\myD_{\myj_1}\cdots\myD_{\myj_{\myk-1}}\myF(\myM^{-1}\myVarX)
	=\sum_{1\leq\myi_1,\ldots,\myi_{\myk-1}\leq\mydim}
	\myD_{\myi_1}\cdots\myD_{\myi_{\myk-1}}\myF(\myVarX)\myME
	\prod_{\myl=1}^{\myk-1}\myME_{\myi_\myl\myj_\myl}.
\end{equation}
Then, applying $D_{\myj_\myk}$ and using
(\ref{lemma::LinearChange::Extended}) again, we obtain
\begin{eqnarray}
\nonumber
	\myD_{\myj_1}\cdots\myD_{\myj_\myk}\myF(\myM^{-1}\myVarX)&=&	\myD_{\myj_\myk}\left(\myD_{\myj_1}\cdots\myD_{\myj_{\myk-1}}\myF\right)(\myM^{-1}\myVarX)
\\
\nonumber
	&=&
	\sum_{\myi_\myk=1}^\mydim\myD_{\myi_\myk}
	\left(\myD_{\myj_1}\cdots\myD_{\myj_{\myk-1}}\myF(\myM^{-1}\myVarX)\right)
	\myME_{\myi_\myk\myj_\myk}
\\
\nonumber
	&=&
	\sum_{\myi_\myk=1}^\mydim\myD_{\myi_\myk}
	\left(
		\sum_{1\leq\myi_1,\ldots,\myi_{\myk-1}\leq\mydim}
		\myD_{\myi_1}\cdots\myD_{\myi_{\myk-1}}\myF(\myVarX)\myME
		\prod_{\myl=1}^{\myk-1}\myME_{\myi_\myl\myj_\myl}
	\right)
	\myME_{\myi_\myk\myj_\myk}
\\
\nonumber
	&=&
	\sum_{1\leq\myi_1,\ldots,\myi_\myk\leq\mydim}
	\myD_{\myi_1}\cdots\myD_{\myi_\myk}\myF(\myVarX)\myME
	\prod_{\myl=1}^\myk\myME_{\myi_\myl\myj_\myl},
\end{eqnarray}
as desired.
\qed

\begin{proposition}\label{prop::DPhiEstimate}
Let $\myk$ and $\nSmooth$, $0\leq\myk\leq\nSmooth$,
be two non negative integers and
let $\myFuncA\in C^0(\R^\mydim)\cap C^\nSmooth(\R^\mydim\setminus\{0\})$
be a $\myMatr-$homogeneous function.
Then, if $1\leq\myj_1,\ldots,\myj_\myk\leq\mydim$,
\begin{eqnarray}
\label{eq::DPhiEstimatesMinus}
\norm{\myVarX}\leq1\ &\Longrightarrow&\ 
\abs{\myD_{\myj_1}\cdots\myD_{\myj_\myk}\myFuncA(\myVarX)}
\leq\bigConstM_1\norm{\myVarX}^{\frac{1}{\mycbase}-\myk}
\\
\label{eq::DPhiEstimatesPlus}
\norm{\myVarX}\geq1\ &\Longrightarrow&\ 
\abs{\myD_{\myj_1}\cdots\myD_{\myj_\myk}\myFuncA(\myVarX)}
\leq\bigConstM_2\norm{\myVarX}^{\frac{1}{\mycmatr}-\myk}
\end{eqnarray}
\end{proposition}

\proof
Let us prove (\ref{eq::DPhiEstimatesMinus}).
Let $\LyapMatrA$ be a \LyapName\ matrix for $\myMatr$
and set
$$
	S_\Lyap=\bigl\{\myVarX\in\R^\mydim\vert\myEucDot[\Lyap\myVarX,\myVarX]=1\bigr\}.
$$
Replacing $\LyapMatrA$ with a $a\LyapMatrA$ for a suitable $a$,
if necessary,
we may assume that for each $\myVarX\in\R^\mydim$
$\myEucDot[\Lyap\myVarX,\myVarX]<1$ if $\norm{\myVarX}\leq1$.

Since $S_\Lyap$ is compact the quantity
$$
	\bigConstM_3=\max\left\{
		\abs{\myD^p\myFuncA(\bar\myVarX)}\mid\bar\myVarX\in S_\Lyap,
		0\leq\abs{p}\leq\nSmooth
	\right\}
$$
is finite.

Let $\myVarX\in\R^\mydim$ and suppose that $\norm{\myVarX}\leq1$.
Choose $\bar\myVarX\in S_\Lyap$ and $\myHomVar\in]0,1]$ such that
$$
	\myVarX=\myHomVar^\myMatr\bar\myVarX.
$$
Denoting by $\myE_1=(1,0,\ldots,0),\ldots,\myE_\mydim=(0,0,\ldots,1)$
the canonical basis of $\R^\mydim$,
from (\ref{prop::EulerEx::ToProve}) it follows that
\begin{eqnarray}
\nonumber
	\myD_{\myj_1}\cdots\myD_{\myj_\myk}\myFuncA(\myVarX)&=&
	\myD_{\myj_1}\cdots\myD_{\myj_\myk}\myFuncA(\myHomVar^\myMatr\bar\myVarX)
\\
\nonumber
	&=&\sum_{1\leq\myi_1,\ldots,\myi_\myk\leq\mydim}
	\myD_{\myi_1}\cdots\myD_{\myi_\myk}\myFuncA(\bar\myVarX)\myHomVar
	\prod_{\myl=1}^\myk\myEucDot[\myE_{\myj_\myl},\myHomVar^{-\myMatr}\myE_{\myi_\myl}],
\end{eqnarray}
and hence, by (\ref{eq::PBaseEstimatesTwo}) and (\ref{eq::PolarB}),
\begin{equation}
\nonumber
	\abs{\myD_{\myj_1}\cdots\myD_{\myj_\myk}\myFuncA(\myVarX)}\leq
	\bigConstM_4\,\myHomVar\prod_{\myl=1}^\myk
	\norm{\myHomVar^{-\myMatr}\myE_\myl}\leq
	\bigConstM_5\,\myHomVar^{1-\mycbase\myk}\leq
	\bigConstM_6\,\norm{\myVarX}^{\frac{1}{\mycbase}(1-\mycbase\myk)}=
	\bigConstM_6\,\norm{\myVarX}^{\frac{1}{\mycbase}-\myk}.
\end{equation}
This proves (\ref{eq::DPhiEstimatesMinus});
an analogous argumentation establishes
(\ref{eq::DPhiEstimatesPlus}).
\qed

In the same way we still have:

\begin{proposition}\label{prop::AHomoEstimates}
Let $\myFuncA\in C^0(\R^\mydim)$ be a positive
$\myMatr-$homogeneous function.
Then, for each $\myVarX\in\R^\mydim$,
\begin{eqnarray}
\norm{\myVarX}\leq1\ &\Longrightarrow&\ 
\bigConstM_1\norm{\myVarX}^{\frac{1}{\mycmatr}}
\leq\myFuncA(\myVarX)\leq
\bigConstM_2\norm{\myVarX}^{\frac{1}{\mycbase}}
\nonumber\\
\norm{\myVarX}\geq1\ &\Longrightarrow&\ 
\bigConstM_3\norm{\myVarX}^{\frac{1}{\mycbase}}
\leq\myFuncA(\myVarX)\leq
\bigConstM_4\norm{\myVarX}^{\frac{1}{\mycmatr}}
\nonumber
\end{eqnarray}
\end{proposition}

\begin{corollary}
Let $\nSmooth>0$ be a positive integer and let
$\myExpA>0$ be a positive real number.
Let $\myFuncA\in C^0(\R^\mydim)\cap C^\nSmooth(\R^\mydim\setminus\{0\})$
be a positive $\myMatr-$homogeneous function.
Assume that
\begin{equation}\nonumber
\myExpA>\mycbase\nSmooth.
\end{equation}
Then
\begin{equation}\nonumber
\myFuncA^\myExpA\in C^\nSmooth(\R^\mydim)
\end{equation}
and for each $\myMultiIndexA\in\N^\mydim$ such that
$\abs\myMultiIndexA\leq\nSmooth$
\begin{equation}\nonumber
\myD^\myMultiIndexA(\myFuncA^\myExpA)(0)=0.
\end{equation}
\end{corollary}

\proof
Observe that
$\myFuncA^\myExpA\in C^0(\R^\mydim)\cap C^\nSmooth(\R^\mydim\setminus\{0\})$
and also it is a $\myExpA^{-1}\myMatr-$homogeneous positive function;
the matrix $\myMatrB=\myExpA^{-1}\myMatr$ satisfies
\begin{eqnarray}
	\myHomVar\geq1&\ \Longrightarrow\ &
		\bigConstM_1\myHomVar^\frac{\mycmatr}{\myExpA}\norm{\myVarX}
		\leq\norm{\myHomVar^{\myMatrB\myVarX}}
		\leq\bigConstM_2\myHomVar^\frac{\mycbase}{\myExpA}\norm{\myVarX},\\
	0<\myHomVar\leq1&\ \Longrightarrow\ &
		\bigConstM_1\myHomVar^\frac{\mycbase}{\myExpA}\norm{\myVarX}
		\leq\norm{\myHomVar^{\myMatrB\myVarX}}
		\leq\bigConstM_2\myHomVar^\frac{\mycmatr}{\myExpA}\norm{\myVarX}.
\end{eqnarray}
Proposition \ref{prop::DPhiEstimate} implies that
if $1\leq\myj_1,\ldots,\myj_\myk\leq\mydim$, $\myk\leq\nSmooth$,
then
\begin{eqnarray}
\norm{\myVarX}\leq1\ &\Longrightarrow&\ 
\abs{\myD_{\myj_1}\cdots\myD_{\myj_\myk}(\myFuncA^\myExpA)(\myVarX)}
\leq\bigConstM_3\norm{\myVarX}^{\frac{\myExpA}{\mycbase}-\myk}
\nonumber
\end{eqnarray}
and hence, being by hypotesis,
\begin{equation}\nonumber
\frac{\myExpA}{\mycbase}-\myk\geq
\frac{\myExpA}{\mycbase}-\nSmooth>0,
\end{equation}
we obtain
\begin{equation}\nonumber
\lim_{x\to0}\myD_{\myj_1}\cdots\myD_{\myj_\myk}(\myFuncA^\myExpA)(\myVarX)=0.
\end{equation}
Since $1\leq\myj_1,\ldots,\myj_\myk\leq\mydim$ and  $\myk\leq\nSmooth$
are arbitrary the assertion easily follows.
\qed

\begin{theorem}\label{thm::ThetaDecay}
Let $\nSmooth>0$ be a positive integer and let
$\myExpA>0$ and $\mydecaya>0$ be positive real numbers.
Let $\myFuncA\in C^0(\R^\mydim)\cap C^\nSmooth(\R^\mydim\setminus\{0\})$
be a positive $\myMatr-$homogeneous function.
Assume that
\begin{equation}\nonumber
\myExpA>\mycbase\nSmooth.
\end{equation}
Then
\begin{equation}\nonumber
e^{-\myFuncA^\myExpA},
\myFuncA^\myExpA e^{-\myFuncA}
\in\mySymSobSpace[\mydecaya,\nSmooth,\mydim].
\end{equation}
\end{theorem}

\proof {
\def\myS{s}
The functions $e^{-\myFuncA^\myExpA}$
and $\myFuncA^\myExpA e^{-\myFuncA}$ are
obviously continuous on $\R^\mydim$.

By the estimates of proposition \ref{prop::AHomoEstimates}
we easily obtain that for $\myVarX\in\R^\mydim$
and $\norm{\myVarX}\geq1$,
\begin{eqnarray}
	&&\abs{e^{-\myFuncA(\myVarX)^\myExpA}}
		\leq e^{-\bigConstM_1\norm{\myVarX}^{\frac{\myExpA}{\mycbase}}},
\nonumber
\end{eqnarray}
and
\begin{eqnarray}
	&&\abs{\myFuncA(\myVarX)^\myExpA e^{-\myFuncA(\myVarX)}}
		\leq\bigConstM_2\norm{\myVarX}^{\frac{\myExpA}{\mycmatr}}
			e^{-\bigConstM_3\norm{\myVarX}^{\frac{1}{\mycbase}}}.
\nonumber
\end{eqnarray}

Since for each $a>0$, $M>0$
\begin{equation}\nonumber
	\lim_{t\to+\infty}t^ae^{-Mt}=0
\end{equation}
it follows that for each $\mydecaya>0$
\begin{equation}\nonumber
e^{-\myFuncA^\myExpA},
\myFuncA^\myExpA e^{-\myFuncA}
\in\myDecaySpace[\mydecaya,\mydim].
\end{equation}

To complete the proof it suffices to prove that
the functions $e^{-\myFuncA^\myExpA}$
and $\myFuncA^\myExpA e^{-\myFuncA}$ are
of class $C^\nSmooth$ and all their derivatives
of order $\myk\leq\nSmooth$ are absolutely integrable
on $\R^\mydim$.

Set $\myFuncB=\myFuncA^\myExpA$.

Let $\myk\leq\nSmooth$.
By induction on $\myk$ it easy to prove that
the derivatives of the function 
$e^{-\myFuncA^\myExpA}=e^{-\myFuncB}$
are linear combinations (with real coefficients) of functions
of the form
\begin{equation}\label{eq::proof::formA}
	\myD^{\myMultiIndexA_1}\myFuncB
	\cdots
	\myD^{\myMultiIndexA_\myS}\myFuncB
	e^{-\myFuncB}
\end{equation}
where
\begin{equation}\nonumber
	\abs{\myMultiIndexA_1}+\cdots+\abs{\myMultiIndexA_\myS}=\myk
\end{equation}
and the derivatives of the function 
$\myFuncA^\myExpA e^{-\myFuncA}=\myFuncB e^{-\myFuncA}$
are linear combinations (with real coefficients) of functions
of the form
\begin{equation}\label{eq::proof::formB}
	\myD^{\myMultiIndexB}\myFuncB
	\myD^{\myMultiIndexA_1}\myFuncA
	\cdots
	\myD^{\myMultiIndexA_\myS}\myFuncA
	e^{-\myFuncA}
\end{equation}
where
\begin{equation}\nonumber
	\abs{\myMultiIndexB}+\abs{\myMultiIndexA_1}+\cdots+\abs{\myMultiIndexA_\myS}=\myk.
\end{equation}

It is now easy to show that each function of the form
(\ref{eq::proof::formA}) either or (\ref{eq::proof::formB})
is $O(\norm{\myVarX}^a)$ for some $a>0$ as $\myVarX\to0$
and $O(\norm{\myVarX}^b e^{-M\norm{\myVarX}^c})$ for some $b,c,M>0$
as $\norm{\myVarX}\to+\infty$.

The proof of the theorem is therefore completed.

} \qed

We end this section with some approximation results.

\begin{proposition}\label{thm::HomogeneousApproximationOne}
Let $\myFuncA_1,\myFuncA_2:\R^\mydim\to[0,+\infty]$ be
two $\myMatr-$homogeneous function.
Assume that $\myFuncA_1$ is upper semicontinuous,
$\myFuncA_2$ is lower semicontinuous
and for each $\myVarX\in\R^n\setminus\{0\}$
$\myFuncA_1(\myVarX)<\myFuncA_2(\myVarX)$.

Then there exists a sequence
$\myFuncB_\myNumA\in C^0(\R^\mydim)\cap C^\infty(\R^\mydim\setminus\{0\})$
of positive $\myMatr-$homogeneous functions such that
for each $\myVarX\in\R^\mydim$
\begin{eqnarray}
&&\myFuncB_{1}(\myVarX)\leq\myFuncA_{2}(\myVarX),
\nonumber\\
&&\myFuncB_{\myNumA+1}(\myVarX)\leq\myFuncB_{\myNumA}(\myVarX),
\quad\myNumA=1,2,\ldots
\nonumber
\end{eqnarray}
and
\begin{eqnarray}
\nonumber&&\lim_{\myNumA\to+\infty}\myFuncB_{\myNumA}(\myVarX)=\myFuncA_{1}(\myVarX).
\end{eqnarray}
\end{proposition}

\proof
Let $\LyapMatrA$ be a \LyapName\ matrix for the matrix $\myMatr$ and set
\begin{equation}\nonumber
	\LyapBall=\bigl\{\myVarX\in\R^\mydim
	\vert\myEucDot[\Lyap\myVarX,\myVarX]=1\bigr\}.
\end{equation}
Then $\LyapBall$ is a compact (sub)manifold.
By standard approximation arguments there exists a non increasing sequence of smooth
positive functions $f_\myNumA$ on $\LyapBall$ which converges pointwise
to the restriction of the function $\myFuncA_{1}$ to $\LyapBall$
and $f_1(\myVarX)\leq\myFuncA_2(\myVarX)$
for each $\myVarX\in\LyapBall$.
For each $\myNumA>0$ let $\myFuncB_\myNumA:\R^\mydim\to[0,+\infty[$
be the unique $\myMatr$-homegeneous functions which extends the function $f_\myNumA$.
Then the sequence $\myFuncB_\myNumA$ has the required properties.

\qed

A similar argument yields:

\begin{proposition}\label{thm::HomogeneousApproximationTwo}
Let $\myFuncA:\R^\mydim\to[0,+\infty]$ be
a continuous $\myMatr-$homogeneous function
and let $0<\e<1$.
Then there exists two positive $\myMatr-$homogeneous functions
$\myFuncB_1,\myFuncB_2\in C^0(\R^\mydim)\cap C^\infty(\R^\mydim\setminus\{0\})$
such that
for each $\myVarX\in\R^\mydim$
\begin{equation}\nonumber
(1-\e)\myFuncA(\myVarX)
\leq\myFuncB_1(\myVarX)
\leq\myFuncA(\myVarX)
\leq\myFuncB_2(\myVarX)
\leq(1+\e)\myFuncA(\myVarX).
\end{equation}
\end{proposition}

} 
\section{\label{section:ProofZetadef}Proof of Theorem \ref{thm::Zetadef}}
	{ 

\def\myHomVarB{u}
\def\mySuperLevel{E}
\def\myRadius{r}
\def\myAuxA{G}


We begin with the following lemma.

\label{eq::ZetaSeries}
	
\begin{lemma}\label{lemma::ZWeakConvergence}
The serie (\ref{eq::ZetaSeries}) defining $\myZeta(\myFuncA,\myComplexVar)$
converges absolutely when $\Re\myComplexVar>\mycbase\mydim$.
\end{lemma}

\proof
If $\myIntVar\in\Z^\mydim\setminus\{0\}$ then $\abs{\myIntVar}\geq1$, and hence,
by Proposition (\ref{prop::AHomoEstimates})
\begin{equation}
\nonumber
	\abs{\myFuncA(\myIntVar)}\geq\bigConstM_1\norm{\myIntVar}^{\frac{1}{\mycbase}}.
\end{equation}
If $\myComplexVar\in\C$ then
\begin{eqnarray}
\nonumber
	&&\abs{\myZeta(\myFuncA,\myComplexVar)}\leq
	\sum_{\myIntVar\in\Z^\mydim\setminus\{0\}}\myFuncA(\myIntVar)^{-\Re\myComplexVar}\leq
	\bigConstM_1^{-1}\sum_{\myIntVar\in\Z^\mydim\setminus\{0\}}
		\norm{\myIntVar}^{-\frac{\Re\myComplexVar}{\mycbase}}
\end{eqnarray}
and the latter series converges (absolutely) when
$\frac{\Re\myComplexVar}{\mycbase}>\mydim$,
that is, if $\Re\myComplexVar>\mycbase\mydim$.
\qed

\begin{proposition}\label{prop::abcZeta}
Let $\myFuncA:\R^\mydim\to[0,+\infty[$ be a continuouos positive
$\myMatr-$homogeneous function.
Let $\myA>0$, $\myB>0$ and $\myC$ be real constants.
Then
\begin{eqnarray}
\label{eq::abcZeta::A}
	\mygbell=e^{-\myFuncA},\ \Re\myComplexVar>\mycbase\mydim
	\ &\Longrightarrow&\ 
	\myXi[\myMatr,\mygbell,\myComplexVar]=
		\Gamma(\myComplexVar)
		\myZeta(\myFuncA,\myComplexVar),
\\
\label{eq::abcZeta::B}
	\mygbell=e^{-\myA\myFuncA},\ \Re\myComplexVar>\mycbase\mydim
	\ &\Longrightarrow&\ 
	\myXi[\myMatr,\mygbell,\myComplexVar]=
		\myA^{-\myComplexVar}
		\Gamma(\myComplexVar)
		\myZeta(\myFuncA,\myComplexVar),
\\
\label{eq::abcZeta::C}
	\mygbell=e^{-\myFuncA^\myB},\ \Re\myComplexVar>\mycbase\mydim
	\ &\Longrightarrow&\ 
	\myXi[\myB^{-1}\myMatr,\mygbell,\myB^{-1}\myComplexVar]=
		\Gamma\left(\frac{\myComplexVar}{\myB}\right)
		\myZeta(\myFuncA,\myComplexVar),
\\
\label{eq::abcZeta::D}
	\mygbell=\myFuncA^{\myC}e^{-\myFuncA},\ \Re\myComplexVar>\max\{\mycbase\mydim,-\myC\},
	\ &\Longrightarrow&\ 
	\myXi[\myMatr,\mygbell,\myComplexVar]=
		\Gamma(\myComplexVar+\myC)
		\myZeta(\myFuncA,\myComplexVar).
\end{eqnarray}
\end{proposition}

\proof
Let $\myC\in\R$ and set $\mygbell=\myFuncA^{\myC}e^{-\myFuncA}$
for $\myVarX\neq0$.
If $\myIntVar\in\Z^\mydim\setminus\{0\}$
and $\Re\myComplexVar>-\myC$ then
\begin{eqnarray}
\nonumber
\int_0^{+\infty}\mygbell(\myHomVar^\myMatr\myIntVar)
	\myHomVar^\myComplexVar\frac{d\myHomVar}{\myHomVar}&=&
	\int_0^{+\infty}\bigl(\myHomVar\myFuncA(\myIntVar)\bigr)^\myC
	e^{-\myHomVar\myFuncA(\myIntVar)}\myHomVar^\myComplexVar\frac{d\myHomVar}{\myHomVar}
\\
\nonumber
&=&\myFuncA(\myIntVar)^\myC
	\int_0^{+\infty}e^{-\myHomVar\myFuncA(\myIntVar)}
	\myHomVar^{\myC+\myComplexVar}\frac{d\myHomVar}{\myHomVar}
\\
\nonumber
&=&\myFuncA(\myIntVar)^\myC
	\int_0^{+\infty}e^{-\myHomVarB}
	\left(\frac{\myHomVarB}{\myFuncA(\myIntVar)}\right)^{\myC+\myComplexVar}\frac{d\myHomVar}{\myHomVar}
\\
\nonumber
&=&\myFuncA(\myIntVar)^{-\myComplexVar}\Gamma(\myC+\myComplexVar).
\end{eqnarray}

By Lemma \ref{lemma::ZWeakConvergence}, if
$\Re\myComplexVar>\max\{\mycbase\mydim,-\myC\}$ then,
summing on $\myIntVar\in\Z^\mydim\setminus\{0\}$,
we obtain
\begin{equation}\nonumber
\myXi[\myMatr,\mygbell,\myComplexVar]=
		\Gamma(\myComplexVar+\myC)
		\myZeta(\myFuncA,\myComplexVar)
\end{equation}
and this proves (\ref{eq::abcZeta::D}).

Setting $\myC=0$ in (\ref{eq::abcZeta::D})
we obtain (\ref{eq::abcZeta::A}).

Let $\myA>0$. Then the assertion (\ref{eq::abcZeta::B}) follows from
(\ref{eq::abcZeta::A}) applied to the $\myMatr-$homogeneous function
$\myFuncB=\myA\myFuncA$, observing that when
$\Re\myComplexVar>\mycbase\mydim$ we trivially have
\begin{equation}\nonumber
	\myZeta(\myFuncB,\myComplexVar)=\myA^{-\myComplexVar}
	\myZeta(\myFuncA,\myComplexVar).
\end{equation}

Let now $\myB>0$.
Then $\myFuncA^\myB$ is a positive $\myB^{-1}\myMatr-$homogeneous function
and the assertion (\ref{eq::abcZeta::C}) follows replacing
in (\ref{eq::abcZeta::A}) the matrix $\myMatr$ with $\myB^{-1}\myMatr$,
$\myFuncA$ with $\myFuncA^\myB$,
and $\myComplexVar$ with $\frac{\myComplexVar}{\myB}$,
observing that when
$\Re\myComplexVar>\mycbase\mydim$ we have
\begin{equation}\nonumber
	\myZeta(\myFuncA^\myB,\frac{\myComplexVar}{\myB})=
	\myZeta(\myFuncA,\myComplexVar).
\end{equation}
\qed

\begin{proposition}\label{prop::GbellIntegral}
Let $\myFuncA\in C^0(\R^\mydim)$ be a positive
$\myMatr-$homogeneous function.
If $\myA>0$ is a positive constant then
\begin{equation}\label{eq::GbellIntegral}
\int_{\R^\mydim}e^{-a\myFuncA(\myVarX)}d\myVarX
=a^{-\mytrace}\Gamma(\mytrace+1)|\myBall_\myFuncA|
\end{equation}
\end{proposition}

\proof
Let $\myVarX\in\R^\mydim$ and let $\myRadius>0$.
By $\myMatr-$homogeneity we have
\begin{equation}\nonumber
	\myFuncA(\myVarX)<\myRadius\ \Longleftrightarrow\ 
	\myRadius^{-1}\myFuncA(\myVarX)<1\ \Longleftrightarrow\ 
	\myFuncA(\myRadius^{-\myMatr}\myVarX)<1,
\end{equation}
that is
\begin{equation}\nonumber
	\myVarX\in\myBall_\myFuncA(\myRadius)\ \Longleftrightarrow\ 
	\myRadius^{-\myMatr}\myVarX\in\myBall_\myFuncA
\end{equation}
and hence
\begin{equation}\label{eq::BallMeasure}
	\abs{\myBall_\myFuncA(\myRadius)}=
	\myRadius^\mytrace\abs{\myBall_\myFuncA}.
\end{equation}

Given $\myHomVar>0$ we set
\begin{equation}\nonumber
\mySuperLevel_\myHomVar=\bigl\{
	\myVarX\in\R^\mydim\mid e^{-a\myFuncA(\myVarX)}>\myHomVar
\bigr\}.
\end{equation}

If $0<\myHomVar<1$ then
\begin{equation}\nonumber
\mySuperLevel_\myHomVar=\myBall_\myFuncA
	\left(
		\frac{1}{\myA}\log\frac{1}{t}
	\right)
\end{equation}
and so, using (\ref{eq::BallMeasure}),
\begin{eqnarray}
	\int_{\R^\mydim}e^{-a\myFuncA(\myVarX)}d\myVarX&=&
	\int_0^1\abs{\mySuperLevel_\myHomVar}d\myHomVar
	=\int_0^1\abs{\myBall_\myFuncA
	\left(
		\frac{1}{\myA}\log\frac{1}{t}
	\right)
	}d\myHomVar
\nonumber\\
	&=&\myA^{-\mytrace}|\myBall_\myFuncA|
	\int_0^1
	\left(
		\log\frac{1}{t}
	\right)^\mytrace
	d\myHomVar
\nonumber\\
	&=&a^{-\mytrace}\Gamma(\mytrace+1)|\myBall_\myFuncA|.
\nonumber
\end{eqnarray}
\qed

Theorem \ref{thm::Zetadef} will be now
an immediate consequence of the following:

\begin{theorem}
Let $\nSmooth$ be a positive integer and assume that $\nSmooth>\mydim$.
Let $\myFuncA\in C^0(\R^\mydim)\cap C^\nSmooth(\R^\mydim\setminus\{0\})$
be a positive $\myMatr-$homogeneous function.
Then the series on the right hand side of (\ref{eq::ZetaSeries})
converges 
to a holomorphic function on the half space
$\Re\myComplexVar>\mytrace$ and extends to a meromorphic function
on the half space $\Re\myComplexVar>\mytrace-\mycmatr\nSmooth$
having only a simple pole at
$\myComplexVar=\mytrace$ with residue
$\mytrace|\myBall_\myFuncA|$.
If $\mytrace-\mycmatr\nSmooth<0$ then we also have
$\myZeta(\myFuncA,0)=-1$.
\end{theorem}

\proof
Let $\myExpA>0$ be a positive real constant such that
$\myExpA>\mycbase\nSmooth$.

Choose $\mydecaya$ satisfying
$\mydecaya>\mydim$,
$\mydecaya>\mycbase\mydim$ and
$\mydecaya+\nSmooth>\frac{\mytrace}{\mycmatr}$.


If $\mygbell=\myFuncA^{\myExpA}e^{-\myFuncA}$ then,
by (\ref{eq::abcZeta::D})
\begin{equation}\nonumber
\mycbase\mydim<\Re\myComplexVar<\mydecaya
	\ \Longrightarrow\
	\myXi[\myMatr,\mygbell,\myComplexVar]=
		\Gamma(\myComplexVar+\myExpA)
		\myZeta(\myFuncA,\myComplexVar).
\end{equation}

By Theorem \ref{thm::ThetaDecay} we also have
$\mygbell\in\mySymSobSpace[\mydecaya,\nSmooth,\mydim]$.

Since $\mygbell(0)=0$,
by (\ref{eq::IntegralEquation}) of Theorem \ref{thm::FunctionalEquation},
when $\mytrace-\mycmatr\nSmooth<\Re\myComplexVar<\mydecaya$
\begin{equation}\nonumber
\myXi[\myMatr,\mygbell,\myComplexVar]=
-\frac{\hat\mygbell(0)}{\mytrace-\myComplexVar}
+\myXiPlus[\myMatr,\mygbell,\myComplexVar]
+\myXiPlus[\myTranspose\myMatr,\hat\mygbell,\mytrace-\myComplexVar],
\end{equation}
the functions 
$\myXiPlus[\myMatr,\mygbell,\myComplexVar]$
and
$\myXiPlus[\myTranspose\myMatr,\hat\mygbell,\mytrace-\myComplexVar]$
being holomorphic when
$\mytrace-\mycmatr\nSmooth<\Re\myComplexVar<\mydecaya$.

Since $\myExpA>\mycbase\mydim$ and $\mycmatr\leq\mycbase$ 
the function $\Gamma(\myComplexVar+\myExpA)$ is holomorphic and not zero
on the strip $\mytrace-\mycmatr\nSmooth<\Re\myComplexVar<\mydecaya$.

It follows that the function
$\Gamma(\myComplexVar+\myExpA)^{-1}\myXi[\myMatr,\mygbell,\myComplexVar]$
is meromorphic on the strip
$\mytrace-\mycmatr\nSmooth<\Re\myComplexVar<\mydecaya$
having only a simple pole at $\myComplexVar=\mytrace$
and it coincide with $\myZeta(\myFuncA,\myComplexVar)$ when
$\mycbase\mydim<\Re\myComplexVar<\mydecaya$;
this shows that the function $\myZeta(\myFuncA,\myComplexVar)$
has a meromorphic extension to the half plane
$\Re\myComplexVar>\mytrace-\mycmatr\nSmooth$
which is holomorphic when $\myComplexVar\neq\mytrace$ and
has a simple pole at $\myComplexVar=\mytrace$
with residue
\begin{equation}\nonumber
	R=\Gamma(\mytrace+\myExpA)^{-1}\hat\mygbell(0)=
		\Gamma(\mytrace+\myExpA)^{-1}
		\int_{\R^\mydim}\myFuncA(\myVarX)^\myExpA e^{-\myFuncA(\myVarX)}d\myVarX.
\end{equation}
Considering $\myExpA$ as a complex variable,
we observe that the function
\begin{equation}\nonumber
	\myExpA\mapsto\myAuxA(\myExpA)=
		\Gamma(\mytrace+\myExpA)^{-1}
		\int_{\R^\mydim}\myFuncA(\myVarX)^\myExpA e^{-\myFuncA(\myVarX)}d\myVarX
\end{equation}
is holomorphic with respect to $\myExpA$ when $\Re\myExpA>0$.
Since $\myAuxA(\myExpA)=R$ when $\myExpA$ is real and $\myExpA>\mycbase\mydim$
then, by the identity principle,
$\myAuxA(\myExpA)=R$ when $\Re\myExpA>0$.

Using the Lebesgue theorem on dominated convergence and (\ref{eq::GbellIntegral})
we obtain
\begin{equation}\nonumber
	R=\lim_{\myExpA\to0^+}\myAuxA(\myExpA)=
	\Gamma(\mytrace)^{-1}\int_{\R^\mydim}e^{-\myFuncA(\myVarX)}d\myVarX
	=\Gamma(\mytrace)^{-1}\Gamma(\mytrace+1)|\myBall_\myFuncA|.
\end{equation}
Recalling that the functional equation for the Euler Gamma function is
\begin{equation}\nonumber
	\Gamma(\mytrace+1)=\mytrace\Gamma(\mytrace),
\end{equation}
it follows that
\begin{equation}\nonumber
	R=\mytrace|\myBall_\myFuncA|.
\end{equation}

Observe now that the series (\ref{eq::ZetaSeries}) defining
$\myZeta(\myFuncA,\myComplexVar)$
converges when $\Re\myComplexVar>\mycbase\mydim$
and extends to a holomorphic function on the bigger half plane
$\Re\myComplexVar>\mytrace$.
The comvergence of the series (\ref{eq::ZetaSeries}) when
$\Re\myComplexVar>\mytrace$
follows then from Landau's Theorem \ref{thm::Landau}.

Assume now that $\mytrace-\mycmatr\nSmooth<0$.
If we set $\mygbell_\myExpA=e^{\myFuncA^\myExpA}$
then, by Theorem \ref{thm::ThetaDecay} we have
$\mygbell_\myExpA\in\mySymSobSpace[\mydecaya,\nSmooth,\mydim]$,
so that by (\ref{eq::abcZeta::C}) and analytic continuation,
if $\mytrace-\mycmatr\nSmooth<\Re\myComplexVar<\mydecaya$,
\begin{equation}\nonumber
	\myZeta(\myFuncA,\myComplexVar)=
	\Gamma\left(\frac{\myComplexVar}{\myExpA}\right)^{-1}
	\myXi[\myExpA^{-1}\myMatr,\mygbell,\myExpA^{-1}\myComplexVar].
\end{equation}
Again by (\ref{eq::IntegralEquation}) of Theorem \ref{thm::FunctionalEquation},
since $\mygbell_\myExpA(0)=1$, in a neighbourhood of $\myComplexVar=0$ we have
\begin{equation}\nonumber
	\myXi[\myExpA^{-1}\myMatr,\mygbell,\myExpA^{-1}\myComplexVar]=
	-\frac{1}{\myComplexVar}+O(1).
\end{equation}
Since
\begin{equation}\nonumber
	\lim_{\myComplexVar\to0}\myComplexVar\Gamma(\myComplexVar)=1
\end{equation}
it follows that
\begin{equation}\nonumber
	\myZeta(\myFuncA,0)=\lim_{\myComplexVar\to0}
	\Gamma\left(\frac{\myComplexVar}{\myExpA}\right)^{-1}
	\myXi[\myExpA^{-1}\myMatr,\mygbell,\myExpA^{-1}\myComplexVar]=
	\lim_{\myComplexVar\to0}
	\Gamma\left(\frac{\myComplexVar}{\myExpA}\right)^{-1}
	\left(
		-\frac{\myExpA}{\myComplexVar}+O(1)
	\right)=-1.
\end{equation}
\qed

} 

The previous Theorem,
combined with Theorem \ref{thm::ThetaDecay} and
Theorem \ref{thm::FunctionalEquation},
immediately yields the following refinement
of Proposition \ref{prop::abcZeta}:

\begin{theorem}\label{thm::abcZeta}
Let $\nSmooth$ be a positive integer and assume that $\nSmooth>\mydim$.
Let $\myFuncA\in C^0(\R^\mydim)\cap C^\nSmooth(\R^\mydim\setminus\{0\})$
be a positive $\myMatr-$homogeneous function.
Let $\myB$ and $\myC$ be positive real constants.
If $\ \Re\myComplexVar>\mytrace-\nSmooth\mycmatr,\ \myComplexVar\neq0,\mytrace,$
then
\begin{eqnarray}
\label{eq::abcZetaEx::C}
	\mygbell=e^{-\myFuncA^\myB},
	\ \myB>\nSmooth\mycbase,
	\ &\Longrightarrow&\ 
	\myXi[\myB^{-1}\myMatr,\mygbell,\myB^{-1}\myComplexVar]=
		\Gamma\left(\frac{\myComplexVar}{\myB}\right)
		\myZeta(\myFuncA,\myComplexVar);
\\
\label{eq::abcZetaEx::D}
	\mygbell=\myFuncA^{\myC}e^{-\myFuncA},
	\ \myC>\nSmooth\mycbase,
	\ &\Longrightarrow&\ 
	\myXi[\myMatr,\mygbell,\myComplexVar]=
		\Gamma(\myComplexVar+\myC)
		\myZeta(\myFuncA,\myComplexVar).
\end{eqnarray}
\end{theorem}

\section{\label{section:ProofZetalim}Proof of Theorem \ref{thm::Zetalim}}
Let $\myMatr$ and $\myFuncA$ be as in the hypotesis of
Theorem \ref{thm::Zetalim}.

When $\varRadius>0$ set
\begin{equation}\nonumber
F_\myFuncA(\varRadius)=\#(\myBall_{\myFuncA_1}(\varRadius)\cap\Z^\mydim)
\end{equation}
and let $\myRadonA_\myFuncA$ be the unique Radon measure 
such that for each $a<b\in\R$
\begin{equation}\nonumber
\mu_\myFuncA(]a,b[)=\lim_{\varRadius\to a^+}F_\myFuncA(\varRadius)-
\lim_{\varRadius\to b^-}F_\myFuncA(\varRadius).
\end{equation}

When $\Re\myComplexVar>\mycbase\mydim$ we have
\begin{equation}\nonumber
\myZeta(\myFuncA,\myComplexVar)=
\int_0^{+\infty}\myHomVar^{-\myComplexVar}d\mu_\myFuncA(\myHomVar).
\end{equation}

Assume that $\myFuncA\in C^0(\R^\mydim)\cap C^\infty(\R^\mydim\setminus\{0\})$
Then Theorem \ref{thm::Zetalim} follows immediately from
Theorem \ref{thm::Zetadef} and the
Ikehara-Wiener Theorem \ref{thm::Ikehara}.

Assume now that $\myFuncA\in C^0(\R^\mydim)$.
Then by theorem \ref{thm::HomogeneousApproximationTwo}
there exist a positive $\myMatr-$homogeneous function
$\myFuncB\in C^0(\R^\mydim)\cap C^\infty(\R^\mydim\setminus\{0\})$
such that for every $\myVarX\in\R^\mydim$
\begin{equation}\nonumber
\myFuncA(\myVarX)\geq\myFuncB(\myVarX).
\end{equation}

Let $\Re\myComplexVar>\mytrace$. Then
\begin{equation}\nonumber
	\abs{\myZeta(\myFuncA,\myComplexVar)}
	\leq\sum_{\myIntVar\in\Z^\mydim\setminus\{0\}}
		\myFuncA(\myIntVar)^{-\Re\myComplexVar}
	\leq\sum_{\myIntVar\in\Z^\mydim\setminus\{0\}}
		\myFuncB(\myIntVar)^{-\Re\myComplexVar}
	<+\infty.
\end{equation}
This complete the proof of assertion 1 of Theorem \ref{thm::Zetalim}.

Fix now $0<\e<1$.
By Theorem \ref{thm::HomogeneousApproximationTwo} there exist
two positive $\myMatr-$homogeneous functions
$\myFuncB_1,\myFuncB_2\in C^0(\R^\mydim)\cap C^\infty(\R^\mydim\setminus\{0\})$
such that for every $\myVarX\in\R^\mydim$
\begin{equation}\nonumber
	(1-\e)\myFuncA(\myVarX)
	\leq\myFuncB_1(\myVarX)
	\leq\myFuncA(\myVarX)
	\leq\myFuncB_2(\myVarX)
	\leq(1+\e)\myFuncA(\myVarX).
\end{equation}

Then we have respectively
\begin{eqnarray}
	&&(1+\e)^{-\mydim}\mya|\myBall_\myFuncA|
	\leq\mya|\myBall_{\myFuncB_2}|
	\leq\liminf_{\sigma\to+\mya^+}(\sigma-\mya)\myZeta({\myFuncB_2},\sigma)
	\nonumber\\
	&&\leq\liminf_{\sigma\to+\mya^+}(\sigma-\mya)\myZeta(\myFuncA,\sigma)
	\leq\limsup_{\sigma\to+\mya^+}(\sigma-\mya)\myZeta(\myFuncA,\sigma)
	\nonumber\\
	&&\leq\limsup_{\sigma\to+\mya^+}(\sigma-\mya)\myZeta({\myFuncB_1},\sigma)
	\nonumber\\
	&&\leq\mya|\myBall_{\myFuncB_1}|
	\leq(1-\e)^{-\mydim}\mya|\myBall_\myFuncA|
	\nonumber
\end{eqnarray}
and
\begin{eqnarray}
	&&(1+\e)^{-\mydim}|\myBall_\myFuncA|\leq|\myBall_{\myFuncB_2}|
	\leq\liminf_{\varRadius\to+\infty}
		\frac{\#(\myBall_{\myFuncB_2}(\varRadius)\cap\Z^\mydim)}{\varRadius^\mya}
	\nonumber\\
	&&\leq\liminf_{\varRadius\to+\infty}
		\frac{\#(\myBall_\myFuncA(\varRadius)\cap\Z^\mydim)}{\varRadius^\mya}
	\leq\limsup_{\varRadius\to+\infty}
		\frac{\#(\myBall_{\myFuncB_1}(\varRadius)\cap\Z^\mydim)}{\varRadius^\mya}
	\nonumber\\
	&&\leq|\myBall_{\myFuncB_1}|
	\leq(1-\e)^{-\mydim}|\myBall_\myFuncA|.
	\nonumber
\end{eqnarray}
Since $\e>0$ can be made arbitrarily small the proof of
assertions 2 and 3 of Theorem \ref{thm::Zetalim} is so completed.

\section{\label{section:ProofZetalim}Proof of Theorem \ref{thm::ZetaDiscontinuous}}
{
\def\mySubSetA{E}
\def\myCompactA{K}
\def\myDenseSet{D}
\def\myDisk{U}
\def\mySurfaceMeas#1{\left|{#1}\right|_{\mydim-1}}

Let $\myMatr$, $\myFuncA$, $\e$ and $\dg$ be as in
Theorem \ref{thm::ZetaDiscontinuous}.

Let $\LyapMatrA$ be a \LyapName\ matrix for the matrix $\myMatr$
and set (again)
\begin{equation}\nonumber
	\LyapBall=\bigl\{\myVarX\in\R^\mydim
	\vert\myEucDot[\Lyap\myVarX,\myVarX]=1\bigr\}.
\end{equation}

For each (Borel) subset $\mySubSetA\sset\LyapBall$ we denote by
$\mySurfaceMeas{\mySubSetA}$ the $(\mydim-1)-$dimensional
Euclidean measure of the set $\mySubSetA$.

Let $\myDenseSet$ be the set of the points
$\myVarX\in\LyapBall$ of the form
$\myHomVar^\myMatr(\myIntVar)$,
where $\myHomVar>0$ and $\myIntVar\in\Z^\mydim$.

As $\myDenseSet$ is countable, it follows that
$\mySurfaceMeas{\myDenseSet}=0$.
By standard measure theory approximation arguments
there exists a compact set
$\myCompactA\sset\LyapBall$ such that
$\myCompactA\cap\myDenseSet=\void$ but
$\mySurfaceMeas{\myCompactA}>0$.

We now define $\myFuncA_1:\R^\mydim\to[0,+\infty[$
as the unique $\myMatr-$homogeneous function
such that when $\myVarX\in\LyapBall$
\begin{equation}\nonumber
	\myFuncA_1(\myVarX)=
	\left\{
		\begin{array}{ll}
			\myFuncA(\myHomVar),&\quad\myVarX\in\LyapBall\setminus\myCompactA,\\
			(1+\e/2)\myFuncA(\myHomVar),&\quad\myVarX\in\myCompactA.
		\end{array}
	\right.
\end{equation}

We also set $\myFuncA_2=(1+\e)\myFuncA$.

Then $\myFuncA_1$ is upper semicontinuous,
$\myFuncA_2$ is lower semicontinuous and
$\myFuncA_1(\myVarX)<\myFuncA_2(\myVarX)$
for each $\myVarX\in\R^\mydim\setminus\{0\}$.

By Proposition \ref{thm::HomogeneousApproximationOne}
there exists a sequence
$$\myFuncB_\myNumA\in C^0(\R^\mydim)\cap C^\infty(\R^\mydim\setminus\{0\})$$ 
of positive $\myMatr-$homogeneous functions such that
for each $\myVarX\in\R^\mydim$,
\begin{eqnarray}
&&\myFuncB_{1}(\myVarX)\leq\myFuncA_{2}(\myVarX)=(1+\e)\myFuncA(\myVarX),
\nonumber\\
&&\myFuncB_{\myNumA+1}(\myVarX)\leq\myFuncB_{\myNumA}(\myVarX),
\quad\myNumA=1,2,\ldots
\nonumber
\end{eqnarray}
and
\begin{eqnarray}
&&\lim_{\myNumA\to+\infty}\myFuncB_{\myNumA}(\myVarX)
=\myFuncA_{1}(\myVarX)\geq\myFuncA(\myVarX).
\nonumber
\end{eqnarray}

It follows that
\begin{eqnarray}
&&\myBall_{\myFuncA_{2}}\sset\myBall_{\myFuncB_{1}}
\nonumber\\
&&\myBall_{\myFuncB_{\myNumA}}\sset\myBall_{\myFuncB_{\myNumA+1}},
\quad\myNumA=1,2,\ldots
\nonumber\\
&&\bigcup_{\myNumA}\myBall_{\myFuncB_{\myNumA}}=\myBall_{{\myFuncA_{1}}}.
\nonumber
\end{eqnarray}
and hence
\begin{eqnarray}
	&&(1+\e)^{-\mydim}|\myBall_\myFuncA|\leq
	|\myBall_{\myFuncB_\myNumA}|\leq|\myBall_\myFuncA|,
	\quad\myNumA=1,2,\ldots
	\nonumber\\
	&&\lim_{\myNumA\to+\infty}|\myBall_{\myFuncB_{\myNumA}}|=
	|\myBall_{\myFuncB_1}|.
	\nonumber
\end{eqnarray}

Since $\mySurfaceMeas{\myCompactA}>0$
and $\mySurfaceMeas{\LyapBall\setminus\myCompactA}>0$
it follows that
\begin{equation}\nonumber
|\myBall_{\myFuncB_1}|<|\myBall_{\myFuncA}|
\end{equation}
and hence
\begin{equation}\label{eq::proof::Inequality}
\lim_{\myNumA\to+\infty}|\myBall_{\myFuncB_{\myNumA}}|
<|\myBall_{\myFuncA}|.
\end{equation}

Observe also that if $\myIntVar\in\Z^\mydim$ then,
by construction,
\begin{eqnarray}\label{eq::proof::IntegerConvergence}
&&\lim_{\myNumA\to+\infty}\myFuncB_{\myNumA}(\myIntVar)
=\myFuncA(\myIntVar).
\end{eqnarray}

Set now
\begin{equation}\nonumber
	\bigConstA_\myNumA=
	\sup_{|\myComplexVar-\mytrace|=\delta}
	\left|\myZeta(\myFuncB_\myNumA,\myComplexVar)\right|.
\end{equation}

We complete the proof of Theorem \ref{thm::ZetaDiscontinuous}
showing that
\begin{equation}\nonumber
	\lim_{\myNumA\to+\infty}\bigConstA_\myNumA=+\infty.
\end{equation}

If not, taking a subsequence if necessary, we have
\begin{equation}\nonumber
	\sup_{|\myComplexVar-\mytrace|=\delta}
	\left|\myZeta(\myFuncB_\myNumA,\myComplexVar)\right|
	\leq\bigConstA
\end{equation}
for some real constant $\bigConstA$.

Consider the sequence of functions defined as
\begin{equation}
	g_\myNumA(\myComplexVar)=\myZeta(\myFuncB_\myNumA,\myComplexVar)
		-\frac{\mytrace|\myBall_{\myFuncB_\myNumA}|}{\myComplexVar-\mytrace}
\end{equation}

By Theorem \ref{thm::Zetadef}
the functions $g_\myNumA(\myComplexVar)$
are holomorphic throughout over all the complex plane.
When $\abs{\myComplexVar-\mytrace}=\delta$ we have
\begin{equation}
	\abs{g_\myNumA(\myComplexVar)}\leq
	\bigConstA+\frac{\mytrace|\myBall_\myFuncA|}{\delta}.
\end{equation}
By the maximun principle the same inequality holds when
$\abs{\myComplexVar-\mytrace}\leq\delta$.

By Vitali-Montel theorem, taking again a subsequence if necessary,
the sequence $g_\myNumA(\myComplexVar)$
converges uniformely on the compact subsets of the disk
$\myDisk=\{\myComplexVar\mid\abs{\myComplexVar-\mytrace}<\delta\}$
to a holomorphic function
$g:\myDisk\to\C$.

Set
\begin{equation}\nonumber
f(\myComplexVar)=g(\myComplexVar)
	+\frac{\mytrace|\myBall_{\myFuncA}|}{\myComplexVar-\mytrace}.
\end{equation}

Since
\begin{equation}\nonumber
\lim_{\myNumA\to+\infty}
	\frac{\mytrace|\myBall_{\myFuncB_\myFuncA}|}{\myComplexVar-\mytrace}
	=\frac{\mytrace|\myBall_{\myFuncA}|}{\myComplexVar-\mytrace}
\end{equation}
uniformely on the compact sets of $\C\setminus\{\mytrace\}$
it follows that 
\begin{equation}\label{eq::proof::ZetaConvergence}
\lim_{\myNumA\to+\infty}
	\myZeta(\myFuncB_\myNumA,\myComplexVar)
	=f(\myComplexVar)
\end{equation}
uniformely on the compact sets of $\myDisk\setminus\{\mytrace\}$.

Assume now that $\myComplexVar=\s\in\R$ and
$\mytrace<\s<\mytrace+\dg$.
Then, by assertion 1 of Theorem \ref{thm::Zetalim},
using (\ref{eq::proof::IntegerConvergence}) and the
Beppo-Levi monotone convergence theorem we obtain
\begin{equation}\nonumber
	f(\s)
	=\lim_{\myNumA\to+\infty}\sum_{\myIntVar\in\Z^\mydim\setminus\{0\}}
			\myFuncB_\myNumA(\myIntVar)^{-\s}
	=\sum_{\myIntVar\in\Z^\mydim\setminus\{0\}}
			\myFuncA(\myIntVar)^{-\s}
	=\myZeta(\myFuncA,\s)
\end{equation}
By assertion 2 of of Theorem \ref{thm::Zetalim}
and using (\ref{eq::proof::ZetaConvergence}) we would obtain
\begin{equation}\nonumber
	\mytrace|\myBall_{\myFuncA}|
	=\lim_{\sigma\to+\mya^+}
		(\sigma-\mya)\myZeta(\myFuncA,\sigma)
	=\myRes_{\myComplexVar=\mytrace}f(\myComplexVar)
	=\lim_{\myNumA\to+\infty}
		\myRes_{\myComplexVar=\mytrace}
			\myZeta(\myFuncB_\myNumA,\myComplexVar)
	=\lim_{\myNumA\to+\infty}\mytrace|\myBall_{\myFuncB_\myNumA}|
\end{equation}
and this contraddicts (\ref{eq::proof::Inequality}).

}
\section{\label{section:ProofThetaAsymp}Asymptotic expantions
and proof of Theorem \ref{thm::ThetaAsymptotic}}
{
\def\myExpA{p}
\def\myM{m}
\def\Region{S}
\def\myBaseB{b}
\def\myBaseT{T}
\def\myRectangleLabel{R}
\def\myRect#1{{\myRectangleLabel(#1)}}
\def\mySector{E}

We begin with an estimate of the growing
of the zeta functions $\myZeta(\myFuncA,\myComplexVar)$
on the imaginary directions.

\begin{proposition}\label{prop::ZetaExpEpsEstimate}
Let $\myFuncA\in C^0(\R^\mydim)\cap C^{\infty}(\R^\mydim\setminus\{0\})$
be a positive $\myMatr-$homogeneous function.

Let $\myA<\myB\in\R$ be given,
and let $\e>0$.
If $\myA\leq\Re\myComplexVar\leq\myB$
and $\abs{\Im\myComplexVar}\geq 1$ then
\begin{eqnarray}
\label{eq::GrowZeta}
	&&\abs{\myZeta(\myFuncA,\myComplexVar)}
	\leq\bigConstM_1e^{\e\Im\myComplexVar},
\\
\label{eq::GrowGammaZeta}
	&&\abs{\Gamma(\myComplexVar)\myZeta(\myFuncA,\myComplexVar)}
	\leq\bigConstM_2e^{-(\pi/2-\e)\Im\myComplexVar}.
\end{eqnarray}
\end{proposition}

\proof
The Stirling formula for the Euler Gamma function implies that
if $\myA\leq\Re\myComplexVar\leq\myB$
and $\abs{\Im\myComplexVar}\geq 1$,
then for each $\dg>0$,
\begin{equation}\label{eq::GrowGamma}
	\bigConstM_3e^{-(\pi/2+\dg)\Im\myComplexVar}
	\leq\abs{\Gamma(\myComplexVar)}
	\leq\bigConstM_4e^{-(\pi/2-\dg)\Im\myComplexVar}
\end{equation}
and hence the estimates (\ref{eq::GrowZeta}) and
(\ref{eq::GrowGammaZeta}) are equivalent.
So it suffices to prove (\ref{eq::GrowGammaZeta}).

Let $\mydecaya>\mydim$ and $\mydecayb>\mydim$
be chosen in such a way
that respectively $\mycmatr\mydecaya>\myB$
and $\mytrace-\mycmatr\mydecayb<\myA$.

Let $\myM$ be any positive integer which satisfies
$\myM>\mycmatr\mydecaya$ wich will be chosed later on.

From Gauss's multiplication formula
\def\myM{m}
\begin{equation}\nonumber
\prod_{\myK=0}^{\myM-1}\Gamma\left(\myComplexVar+\frac{\myK}{\myM}\right)
=\myM^{\frac{1}{2}-\myM\myComplexVar}
(2\pi)^{\frac{1}{2}(\myM-1)}\Gamma(\myM\myComplexVar)
\end{equation}
it follows that
\begin{equation}\nonumber
\Gamma(\myM\myComplexVar)=
(2\pi)^{\frac{1}{2}(1-\myM)}
\myM^{\myM\myComplexVar-\frac{1}{2}}
\Gamma(\myComplexVar)
\prod_{\myK=1}^{\myM-1}\Gamma\left(\myComplexVar+\frac{\myK}{\myM}\right);
\end{equation}
then
\begin{equation}\nonumber
\Gamma(\myComplexVar)=
(2\pi)^{\frac{1}{2}(1-\myM)}
\myM^{\myComplexVar-\frac{1}{2}}
\Gamma\left(\frac{\myComplexVar}{\myM}\right)
\prod_{\myK=1}^{\myM-1}\Gamma\left(\frac{\myComplexVar+\myK}{\myM}\right).
\end{equation}

Setting $\mygbell_\myM(\myVarX)=e^{-\myFuncA(\myVarX)^\myM}$ and
using the identity (\ref{eq::abcZetaEx::C})
of Theorem \ref{thm::abcZeta} we obtain

\begin{equation}\nonumber
\Gamma(\myComplexVar)\myZeta(\myFuncA,\myComplexVar)=
\myXi[\myM^{-1}\myMatr,\mygbell,\frac{\myComplexVar}{\myM}]
(2\pi)^{\frac{1}{2}(1-\myM)}
\myM^{\myComplexVar-\frac{1}{2}}
\prod_{\myK=1}^{\myM-1}\Gamma\left(\frac{\myComplexVar+\myK}{\myM}\right).
\end{equation}

The identity (\ref{eq::IntegralEquation})
of Theorem \ref{thm::FunctionalEquation} yields therefore

\begin{equation}\nonumber
\myXi[\myM^{-1}\myMatr,\mygbell,\frac{\myComplexVar}{\myM}]=
-\frac{\myM}{\myComplexVar}
-\frac{\myM\hat\mygbell_\myM(0)}{\mytrace-\myComplexVar}
+\myXiPlus[\myM^{-1}\myMatr,\mygbell_\myM,\frac{\myComplexVar}{\myM}]
+\myXiPlus[\myM^{-1}\myTranspose\myMatr,\hat\mygbell_\myM,
		\frac{\mytrace-\myComplexVar}{\myM}].
\end{equation}

Assume now that $\myComplexVar$
lies in the region $\Region$
defined by the conditions
$\myA\leq\Re\myComplexVar\leq\myB$
and $\abs{\Im\myComplexVar}\geq 1$.
Then the estimate (\ref{eq::XiPlusEstimate})
of Lemma \ref{prop::XiPlusEstimate} 
implies
\begin{equation}\nonumber
\abs{\myXiPlus[\myM^{-1}\myMatr,\mygbell_\myM,\frac{\myComplexVar}{\myM}]}
\leq
\frac{\myM\bigConstM_5\myDecayNorm[\mygbell_\myM,\mydecaya]}{\mycmatr\mydecaya-\Re\myComplexVar}
\leq
\frac{\myM\bigConstM_5\myDecayNorm[\mygbell_\myM,\mydecaya]}{\mycmatr\mydecaya-\myB}
\end{equation}
and
\begin{equation}\nonumber
\abs{\myXiPlus[\myM^{-1}\myTranspose\myMatr,\hat\mygbell_\myM,
		\frac{\mytrace-\myComplexVar}{\myM}]}
\leq
\frac{\myM\bigConstM_6\myDecayNorm[\mygbell_\myM,\mydecayb]}{\Re\myComplexVar-\mytrace+\mycbase\mydecayb}
\leq
\frac{\myM\bigConstM_6\myDecayNorm[\mygbell_\myM,\mydecayb]}{\myA-\mytrace+\mycbase\mydecayb}.
\end{equation}

It follows easily now that the function
$\myXi[\myM^{-1}\myMatr,\mygbell,\frac{\myComplexVar}{\myM}]$
is bounded on the region $\Region$
and hence
\begin{equation}\nonumber
\abs{\Gamma(\myComplexVar)\myZeta(\myFuncA,\myComplexVar)}
\leq
\bigConstM_7
(2\pi)^{\frac{1}{2}(1-\myM)}
\myM^{\myB-\frac{1}{2}}
\prod_{\myK=1}^{\myM-1}
\abs{\Gamma\left(\frac{\myComplexVar+\myK}{\myM}\right)}
=\bigConstM_8
\prod_{\myK=1}^{\myM-1}
\abs{\Gamma\left(\frac{\myComplexVar+\myK}{\myM}\right)}.
\end{equation}

Setting $\dg=\e/2$ in (\ref{eq::GrowGamma}) we easily obtain
\begin{equation}\nonumber
\abs{\Gamma(\myComplexVar)\myZeta(\myFuncA,\myComplexVar)}
\leq
\bigConstM_9
e^{-\left(\frac{\pi-\e}{2}\right)
	  \left(\frac{\myM-1}{\myM}\right)
	  \Im\myComplexVar}.
\end{equation}
We now chose $\myM$ large enough to satisfy
\begin{equation}\nonumber
\left(\frac{\pi-\e}{2}\right)
	  \left(\frac{\myM-1}{\myM}\right)
>\frac{\pi}{2}-\e
\end{equation}
and we are done.
\qed

\begin{proposition}
Let $\myFuncA\in C^0(\R^\mydim)$
be a positive $\myMatr-$homogeneous function.
If $\myHomVar>0$ is a positive real number
\begin{equation}\nonumber
	\theta(\myFuncA,i\myHomVar)
	=\myTheta[\myMatr,e^{-\myFuncA},i\myHomVar]
\end{equation}
and if $\Im\tau>0$,
\begin{equation}\nonumber
	\abs{\theta(\myFuncA,\tau)}
	\leq\theta(\myFuncA,\Im\tau).
\end{equation}
\end{proposition}

\proof
When $\myHomVar>0$, by $\myMatr-$homogeneity, we have
\begin{equation}\nonumber
	\theta(\myFuncA,i\myHomVar)
	=\sum_{\myIntVar\in\Z^\mydim}e^{-\myHomVar\myFuncA(\myIntVar)}
	=\sum_{\myIntVar\in\Z^\mydim}e^{-\myFuncA(\myHomVar^\myMatr\myIntVar)}
	=\myTheta[\myMatr,e^{-\myFuncA},i\myHomVar].
\end{equation}

When $\Im\tau>0$
\begin{eqnarray}\nonumber
	\abs{\theta(\myFuncA,\tau)}&\leq&
	\sum_{\myIntVar\in\Z^\mydim}\abs{e^{i\tau\myFuncA(\myIntVar)}}
	\leq
	\sum_{\myIntVar\in\Z^\mydim}e^{-\Im\tau\myFuncA(\myIntVar)}
	=\theta(\myFuncA,\Im\tau).
\end{eqnarray}
\qed

We are now ready to prove Theorem \ref{thm::ThetaAsymptotic}.

Let $\Re\rhVar>0$ and
$\Re\myComplexVar>\mycbase\mydim$.
Then
\begin{equation}\nonumber
	\rhVar^{-\myComplexVar}\Gamma(\myComplexVar)\myZeta(\myFuncA,\myComplexVar)=
	\int_0^{+\infty}\theta^*(\myFuncA,i\myHomVar\rhVar)\myHomVar^\myComplexVar
		\frac{d\myHomVar}{\myHomVar}.
\end{equation}

Indeed both sides are holomorphic
with respect to the variable $\rhVar$
on the half plane $\Re\rhVar>0$ 
and coincide when $\rhVar=\myA>0$ 
by (\ref{eq::abcZeta::B}) of proposition \ref{prop::abcZeta}.

Using the estimates for
$\abs{\Gamma(\myComplexVar)\myZeta(\myFuncA,\myComplexVar)}$
given in proposition \ref{prop::ZetaExpEpsEstimate}
we easily obtain that for each $\e>0$
\begin{equation}\label{eq::proof::Grow}
	\abs{\rhVar^{-\myComplexVar}\Gamma(\myComplexVar)\myZeta(\myFuncA,\myComplexVar)}
	=O(e^{-(\pi/2-\abs{\arg{\rhVar}}-\e)})
\end{equation}
as $\abs{\Im\myComplexVar}\to+\infty$ uniformely with respect to $\Re\myComplexVar$.
The Mellin inversion formula implies
therefore that, when $\myBaseB>\mycbase\mydim$,
\begin{equation}\label{eq::proof::Mellin}
	\theta(\myFuncA,i\rhVar)=
	1+\theta^*(\myFuncA,i\rhVar)=
	1+\frac{1}{2\pi i}
	\int_{(\myBaseB)}\rhVar^{-\myComplexVar}\Gamma(\myComplexVar)
		\myZeta(\myFuncA,\myComplexVar)d\myComplexVar.
\end{equation}
Let $\myN$ be a positive integer and let $0<\vLine<1$.
Given a positive real number $\myBaseT>0$ consider the rectangle
$\myRect\myBaseT$ with vertices
$-\myN-1+\vLine\pm i\myBaseT$ and $\myBaseB\pm i\myBaseT$.
The poles of the function
$\rhVar^{-\myComplexVar}\Gamma(\myComplexVar)\myZeta(\myFuncA,\myComplexVar)$
inside $\myRect\myBaseT$ are
$\myComplexVar=\mytrace$ with residue
$
R=\Gamma(\mytrace+1)\abs{\myBall_\myFuncA}\rhVar^{-\mytrace}$,
$\myComplexVar=0$ with residue
$R=\myZeta(\myFuncA,0)=-1$,
and for each integer $\myK=1,\ldots,\myN$,
$\myComplexVar=-\myK$ with residue
$R=((-1)^\myK\myZeta(\myFuncA,-\myK)/\myK!)\rhVar^{\myK}$.

The residue theorem implies that
\begin{equation}\nonumber
	\frac{1}{2\pi i}\int_{\partial\myRect\myBaseT}
		\rhVar^{-\myComplexVar}\Gamma(\myComplexVar)
		\myZeta(\myFuncA,\myComplexVar)d\myComplexVar
	=\Gamma(\mytrace+1)\abs{\myBall_\myFuncA}\rhVar^{-\mytrace}
	-1+
	\sum_{\myK=1}^{\myN}
		\frac{(-1)^\myK\myZeta(\myFuncA,-\myK)}{\myK!}\rhVar^\myK.
\end{equation}

Letting $\myBaseT\to+\infty$,
the estimate (\ref{eq::proof::Grow}) implies that
\begin{eqnarray}\nonumber
	\frac{1}{2\pi i}\int_{(\myBaseB)}
		\rhVar^{-\myComplexVar}\Gamma(\myComplexVar)
		\myZeta(\myFuncA,\myComplexVar)d\myComplexVar
	&=&
	\frac{1}{2\pi i}\int_{(-N-1+\vLine)}
		\rhVar^{-\myComplexVar}\Gamma(\myComplexVar)
		\myZeta(\myFuncA,\myComplexVar)d\myComplexVar
\\
\nonumber
	&+&\Gamma(\mytrace+1)\abs{\myBall_\myFuncA}\rhVar^{-\mytrace}
	-1+
	\sum_{\myK=1}^{\myN}
		\frac{(-1)^\myK\myZeta(\myFuncA,-\myK)}{\myK!}\rhVar^\myK.
\end{eqnarray}

Inserting such expression in (\ref{eq::proof::Mellin})
we obtain (\ref{eq::ThetaAsymptotic}).

The estimate (\ref{eq::ThetaAsymptoticEstimate})
is an easy consequence of
(\ref{eq::proof::Grow})
and the proof is complete.

\rem Classical Example.::
Let $0<\vSector<\pi/2$.
Let denote by $\mySector(\vSector)$ the open angle
defined by the equations
$\Re\rhVar>0$ and $\abs{\arg\rhVar}<\pi/2-\vSector$.
Let $\myExpA>0$.
For $\mydim=1$, $\myFuncA(\myVarX)=\abs{\myVarX}^\myExpA$ we have
$\myZeta(\myFuncA,\myComplexVar)=2\myZeta(\myExpA\myComplexVar)$,
where $\myZeta(\myComplexVar)$ is the Riemann zeta function
and $\abs{\myBall_\myFuncA}=\abs{[-1,1]}=2$;
in this case Theorem \ref{thm::ThetaAsymptotic} yields
\begin{equation}\nonumber
\sum_{\myIntVar=-\infty}^{+\infty}
	e^{-\rhVar\abs{\myIntVar}^\myExpA}=
2\Gamma\left(\frac{1}{\myExpA}+1\right)\rhVar^{-1/\myExpA}+
2\sum_{\myK=1}^{\myN}
	\frac{(-1)^\myK\myZeta(-\myK\myExpA)}{\myK!}\rhVar^\myK+
O(\abs{\rhVar}^{\myN+1-\vLine}),
\quad\rhVar\in\mySector(\vSector).
\end{equation}

Since the Riemann zeta function has its real zeroes
at the even negative integers,
when $\myExpA=2m$ is an even positive integer
we obtain that for each $\mydecaya>0$
\begin{equation}\nonumber
\frac{1}{2}
\sum_{\myIntVar=-\infty}^{+\infty}
	e^{-\rhVar\myIntVar^{2m}}=
\Gamma\left(\frac{1}{2m}+1\right)\rhVar^{-1/{2m}}+
O(\abs{\rhVar}^{\mydecaya}),
\quad\rhVar\in\mySector(\vSector).
\end{equation}

When $\myExpA=1$ we obtain easily
\begin{equation}\nonumber
\frac{1}{2}\,\frac{1+e^{-\rhVar}}{1-e^{-\rhVar}}=
\frac{1}{\rhVar}+
\sum_{\myK=1}^{\myN}
	\frac{(-1)^\myK\myZeta(-\myK)}{\myK!}\rhVar^\myK+
O(\abs{\rhVar}^{\myN+1-\vLine}),
\end{equation}
according to the fact that when $\myK=1,2,\ldots,$ 
$-(\myK+1)\myZeta(-\myK)=B_{\myK+1}$,
where $B_{\myK+1}$ is the $(\myK+1)-$th Bernoulli number.

}


\eject

\nocite{book:ArnoldOrdDiffEqSpringer}
\nocite{book:ArnoldOrdDiffEqMIT}
\nocite{book:LangAlgebraicNumerTheory}
\nocite{book:TenenbaumAnalyticNumberTheory}
\nocite{article:EpsteinA}
\nocite{article:EpsteinB}
\nocite{article:SiegelA}
\nocite{article:SiegelB}

\end{document}